\documentclass[english,letterpaper,11pt,reqno]{amsart}

\usepackage{appendix}
\usepackage{amsbsy}
\usepackage{amsfonts}
\usepackage{amsmath}
\usepackage{amssymb}
\usepackage{amsthm}
\usepackage{graphicx}
\usepackage{ifthen}
\usepackage{textcomp}
\usepackage{enumitem, color, amssymb}
\usepackage{mathrsfs}
\usepackage[bookmarksnumbered,colorlinks]{hyperref}
\hypersetup{colorlinks=true, linkcolor=blue, citecolor=blue}
\usepackage{placeins}

\newtheorem {lemma} {Lemma} [section]
\newtheorem{thm}{Theorem}

\theoremstyle{remark}
\newtheorem{remark}[lemma]{Remark}

\newcommand{\beqa}{\begin{eqnarray}}
\newcommand{\beq}{\begin{equation}}
\newcommand{\eeqa}{\end{eqnarray}}
\newcommand{\eeq}{\end{equation}}

\newcommand{\be}{\begin{equation}}
\newcommand{\ee}{\end{equation}}
\newcommand{\lb}[1]{\label{#1}}

\renewcommand{\Ref}[1]{(\ref{#1})}

\newcommand\kk{{\boldsymbol k}}

\newcommand\vv[1]{{\boldsymbol {\it #1}}}
\newcommand\xx{\vv{x}}
\newcommand\yy{\vv{y}}
\newcommand\zz{\vv{z}}

\newcommand\n{\nabla}
\newcommand\gK{g_{\scriptscriptstyle K}}
\newcommand\nK{\n^{\scriptscriptstyle K}}

\newcommand{\HH}{\mathcal{H}}
\newcommand{\VV}{\mathcal{V}}
\newcommand{\ta}{\tau}
\newcommand\om{\omega}
\newcommand\tT{{\boldsymbol t}}
\newcommand\pp{{\boldsymbol p}}

\newcommand{\al}{\alpha}
\newcommand{\bet}{\beta}
\newcommand{\lam}{\lambda}

\newcommand{\we}{\wedge}

\newcommand{\fr}{\frac}

\begin{document}
\title[]{Canonical K\"ahler metrics on classes of Lorentzian $4$-manifolds}
\author[]{Amir Babak Aazami and Gideon Maschler}
\address{Department of Mathematics and Computer Science\\ Clark University\\ Worcester, MA }
\email{Aaazami@clarku.edu\,,\,Gmaschler@clarku.edu}

\maketitle
\thispagestyle{empty}
\begin{abstract}
Conditions for the existence of K\"ahler-Einstein metrics and central K\"ahler metrics~\cite{cen}
along with examples, both old and new, are given on classes of Lorentzian $4$-manifolds with two distinguished vector fields. The results utilize the general construction \cite{a-m} of K\"ahler metrics on such manifolds. The examples include both complete and incomplete metrics, and some
reside on Lie groups associated to four types of Lie algebras. An appendix includes a similar
construction for scalar-flat K\"ahler metrics.
\end{abstract}

\section{Introduction}
\lb{sec:0}
In \cite{a-m}, the authors introduced a construction scheme associating a family
of K\"ahler metrics to an oriented Lorentzian (or even semi-Riemannian) $4$-manifold
equipped with data tied to two distinguished vector fields. Many examples
were given. Included among the Lorentzian metrics to which the construction applies
are certain warped products such as de Sitter spacetime, gravitational plane waves,
some Petrov type $D$ metrics such as the Kerr metric, metrics which yield the class
of K\"ahler metrics known as SKR metrics (see \cite{dm1}), and a certain
solvable Lie group. The K\"ahler metrics induced by this construction
may, in general, be given only in subregions of the original $4$-manifold, but in many
cases they are, in fact, everywhere defined.

Our purpose in this paper is to study cases where this construction leads to
{\em curvature-distinguished} K\"ahler metrics. In \cite{a-m} one such case was
noted, namely when the K\"ahler metric within the SKR class is the extremal K\"ahler
metric conformal to Page's Einstein metric \cite{p}. Curvature computations in \cite{a-m}
that were relevant to that case were given in part with respect to coordinates.
The methods of the current paper are based on frames rather than coordinates.

The metrics we describe are not all new, as a main part of our motivation is to demonstrate
how familiar metrics arise from the Lorentzian setting described in \cite{a-m}. Specifically, we exhibit
classes of central K\"ahler metrics, K\"ahler-Einstein metrics, and in an appendix,
some scalar-flat K\"ahler surfaces, all on non-compact $4$-manifolds. We now describe these in more detail.

{\em Central K\"ahler metrics}, are metrics for which the determinant of the Ricci endomorphism is constant. On compact manifolds these were studied in \cite{cen} in analogy with constant
scalar curvature K\"ahler metrics, and they share some of the latter's structure, including
an associated Futaki invariant. In the present context we exhibit new metrics of zero Ricci determinant,
which can also be described on 4-manifolds as metrics having a semi-definite Ricci tensor. It is interesting to note that for compact manifolds, zero Ricci determinant metrics turn out to be the hardest to classify.

Theorems~\ref{central} and \ref{KEc} describe central metrics and K\"ahler-Einstein metrics
induced by certain Lorentzian metrics admitting either a Killing field or a vector field with
a geodesic flow. The main Lorentzian metrics we give as examples of this construction
are plane gravitational waves,  metrics related to other pp-waves or metrics residing
on four simply-connected Lie groups associated to the Lie algebras $\mathfrak{su}(2)\times\mathbb{R}$, $\mathfrak{sl}(2)\times\mathbb{R}$, $\mathfrak{nil}_3\times\mathbb{R}$, as well as a fourth,
less well-known step 3 nilpotent Lie algebra.

When a Lorentzian invariant of one of the vector fields in the construction, called the twist function, satisfies a certain differential equation,  we show these central metrics are conformal to a metric of constant scalar curvature (CSC), whereas the K\"ahler-Einstein metrics are only biconformal to such a metric. In the special case where the twist function is constant, this conformally related metric is
locally isometric to a left-invariant metric on a Lie group. In the above mentioned Lie group
examples the CSC metric is in fact left-invariant, hence, in particular, complete.

Another class of K\"ahler-Einstein metrics is described in Theorem~\ref{ke}. These arise from nontrivial
Lorentzian warped products of the real line with a $3$-dimensional fiber satisfying certain conditions.
Theorem~\ref{ke} shows that these metrics arise precisely when both a differential equation holds for a
similar twist function on the fiber, and an additional ODE is satisfied.
These K\"ahler-Einstein metrics belong to a class produced by a classical ansatz given in
\cite{bb}, that has seen many more recent generalizations (see for example \cite{da-wa1}, \cite{da-wa2}, \cite{wa}, \cite{ap-ca-gau}).

The well-known ansatz of LeBrun~\cite{leb} for producing scalar-flat K\"ahler metrics on $4$-manifolds
is discussed in the appendix from a viewpoint closely related to \cite{a-m} and the rest of the current
work. We show how to obtain examples starting from an exact symplectic form associated to a Lorentzian
or other semi-Riemannian metric, and give a few examples where the latter is a pp-wave, and the scalar
flat K\"ahler metrics are also conformal to an Einstein metric.




Finally, we note that our examples include both complete and incomplete metrics. For the case of
the K\"ahler-Einstein metrics of Theorem~\ref{KEc}, our methods allow one to reproduce the complete
biaxial $SU(2)$ cohomogeneity one metrics in \cite{g-p}, \cite{pe}, \cite{d-s1}, \cite{d-s2}.
In the latter reference a local cohomogeneity one metric for the Heisenberg group is given without
results on completeness. We show that a complete cohomogeneity one K\"ahler-Einstein metric exists under
the action of a quotient group of the Heisenberg group. These methods are in part based on recent characterizations of metric smoothness at singular orbits, given in \cite{vz}. We hope to extend these arguments to other groups
in future work.



Section~\ref{sec:1} recalls the construction in \cite{a-m}. In Section~\ref{sec:2}
we first relay the result on central metrics (Theorem~\ref{central}), give examples,
develop the theory on Lie groups and then describe the K\"ahler-Einstein metrics
involved in this construction (Theorem~\ref{KEc}), including completeness.
In Section~\ref{sec:3}  we give Theorem~\ref{ke} concerning our other
class of K\"ahler-Einstein metrics. The appendix describes our result concerning
scalar-flat K\"ahler metrics.

\medskip
\noindent
{\bf Acknowledgements}
We are grateful to Andrzej Derdzinski and Robert Ream for frutiful exchanges regarding
Lie groups and completeness; more details are given in the text.
We also thank the referee for pointed remarks that led to an expansion of an earlier
version to include the material on Lie group metrics, Theorem~\ref{KEc}, completeness and
the appendix, as well as improvements to other aspects of this work, in particular subsection~\ref{k-dep}.


\subsection{$k$-dependency}\lb{k-dep}
This subsection, which is not strictly necessary for following the derivation of
our theorems, outlines a general principle which, in a sense, enables them. It
is described via the following notion of {\em $k$-dependency},
which can be viewed as a local version of the notion of
metric symmetry or cohomogeneity, that is defined via frames.
We first illustrate the principle in a special case.

Given a frame $\{e_i\}$, a function $\ta$ and a semi-Riemannian metric $g$
on a manifold $M$, suppose
\[ d_{e_i}\ta=a_i(\ta),\qquad [e_i,e_j]=b_{ij}^k(\ta)e_k,\qquad g(e_i,e_j)=c_{ij}(\ta), \]
where $d_{e_i}$ denote frame directional derivatives and the summation convention
is employed in the second equality (and again below).
Here $a_i(\ta)=\tilde{a}_i\circ\ta$, with the $\tilde{a}_i$ being
real-valued functions defined on the range of $\ta$, and
similar definitions hold for $b_{ij}^k(\ta)$ and $c_{ij}(\ta)$.
In this setting we will say that $(M,g)$ is $1$-dependent, this
referring to all the above quantities being dependent on {\em one function} $\ta$.
The main utility of this definition stems from the fact that the Koszul formula
involves only Lie brackets, the metric and its directional derivatives.
It follows from this that evaluating a covariant derivative in a frame field's
direction on another frame field, will yield a vector field whose expression as
a combination of the frame fields will have coefficients that are also $\tau$-dependent.
Explicitly, by Koszul's formula,
\begin{multline*}
2g(\n_{e_i}e_j,e_k)=c_{jk}'(\ta)a_i(\ta)+c_{ik}'(\ta)a_j(\ta)-c_{ij}'(\ta)a_k(\ta)\\
-b_{jk}^l(\ta)c_{il}(\ta)-b_{ik}^l(\ta)c_{jl}(\ta)+b_{ij}^l(\ta)c_{kl}(\ta)=:2A_{ij}^k(\ta),
\end{multline*}
where $c_{jk}'(\ta)=\tilde{c}_{jk}'\circ\ta$ etc., with the prime denoting $d/d\ta$.
Therefore, if $\n_{e_i}{e_j}=B_{ij}^\ell e_\ell$, then the $B_{ij}^\ell$ are obtained
by solving the linear system $A_{ij}^k(\ta)=B_{ij}^\ell c_{\ell k}(\ta)$, so they are
also functions of $\ta$.

It is easily seen that the same conclusion  regarding $\ta$-dependence applies to all curvature quantities evaluated on the frame vector fields. For example, for the curvature we have
\begin{multline*}
R(e_i,e_j)=-[{B_{jk}^m}'(\ta)a_i(\ta)-{B_{ik}^m}'(\ta)a_j(\ta)+\\
B_{jk}^\ell(\ta)B_{i\ell}^m(\ta)
-B_{ik}^\ell(\ta)B_{jl}^m(\ta)
-b_{ij}^\ell(\ta)B_{lk}^m(\ta)]e_m.
\end{multline*}
Therefore, natural curvature conditions such
as the Einstein condition, when evaluated on frame fields, will result in a complicated system of ODEs,
where the independent variable is $\ta$ and the dependent variables are the tilded quantities.
A solution of this system gives the information necessary to construct the metric in the given frame
domain. Depending on the problem at hand, one may assume, for example, that $a_i(\ta)$ and $b_{jk}^l(\ta)$ are known manifold data, in which case they become given coefficients in the equations, and attempt to solve for the $c_{ij}(\ta)$.


In analogy with this scheme, given two functionally independent functions $\ta$, $\mu$ on $M$, one can instead require $a_i$, $\bar a_i$, $b_{ij}^k$, $c_{ij}$, with $\bar a_i=d_{e_i}\mu$, to have the form $\tilde{a}_i\circ (\ta,\mu)$, $\tilde{\bar{a}}_i\circ (\ta,\mu)$,
$\tilde{b}_{ij}^k\circ(\ta,\mu)$, $\tilde{c}_{ij}\circ(\ta,\mu)$ respectively, where the tilded quantities are functions on the range of $(\ta,\mu)$ in $\mathbb{R}^2$. This will be the 2-dependent case. Natural curvature conditions applied to the frame fields will now yield PDEs with independent
variables $\ta$, $\mu$ and tilded functions the dependent variables.

In general, we call a semi-Riemannian manifold {\em $k$-dependent} if there exists local frames near each point such that for each frame, there are $k$ functionally independent real-valued functions $\{\ta_i\}$ on the frame domain, called henceforth
the $k$-set, whose frame directional derivatives $a_i^{\ta_j}$ are functionally dependent on the $k$-set,
and so are
the metric values $c_{ij}$ on the frame vector fields 
and the coefficients $b_{ij}^k$ in the expansion via the frame fields of all Lie brackets of pairs of frame fields. In the special case where $c_{ij}$, $b_{ij}^k$ are constant, we call the semi-Riemannian
manifold $0$-dependent.

Once more, in the $k$-dependent case, the Koszul formula shows that the connection, curvature
and Ricci curvature components on the given frame fields also depend functionally on the
$k$-set, and any natural curvature condition, such as the Einstein condition, is equivalent locally
to a PDE system in which the $k$-set functions serve as the independent variables, whereas
depending on what is considered given, the unknowns are taken from the tilded quantities associated to the $a_i^{\ta_j}$'s, $b_{ij}^k$'s and $c_{ij}$'s.

This concept of $k$-dependency is clearly sensible when $k<\dim M$. Two important examples are as follows. In the $0$-dependent case the metric is locally isometric to a left-invariant metric on some
Lie group. The second example is given by a coordinate frame in which $p$ of the coordinates vector
fields are Killing, which gives rise to $(n-p)$-dependency on the coordinate domain, where $n$ is the
manifold dimension. Here the frame is the coordinate frame and the $(n-p)$-set consists of coordinate
functions whose coordinate vector fields are not Killing. In this case, frame directional derivatives of the functions in the $(n-p)$-set are constant, equal to either $1$ or $0$. Additionally, all the Lie brackets of frame fields vanish, whereas the metric coefficients depend only on the coordinate functions belonging to the $(n-p)$-set and not on the others. More generally, one can take instead of such an abelian Lie algebra
of Killing fields, any Lie algebra of Killing fields acting freely.


How the results of this paper actually fit in this framework is not completely straightforward.
In all our theorems the $b_{ij}^k$'s and $c_{ij}$'s depend on two functions,
denoted $\ta$ and either $\iota$ or $\bar{\iota}$, the latter two being twist functions.
Now the  $a_i^{\ta}$'s are constants, so that in the case such a twist function is constant,
we are in the setting of $1$-dependency (dependence on $\ta$), and the relevant theorem translates the
curvature condition into a single ODE (although in Theorem~\ref{central} this ODE does not appear
explicitly, as we effectively express the theorem directly in terms of the ODE's solution).

In the case that the twist function is not constant, we need to examine the $a_i^\iota$ (or the $a_i^{\bar{\iota}}$ for Theorem~\ref{ke}) to determine the dependency type. However, not all
of them are given explicitly, and therefore in that case the curvature condition translates
into an additional ``generalized PDE".

To explain this term, let $\{e_i\}$, $i=1\ldots n$ be an ordered frame on a manifold $M$.
Denote by $d_{e_i}$ a directional derivative with respect to
a frame vector field $e_i$, while its m-fold composition, $m\ge 0$, will be
denoted $d_{e_i}^m$. If $\ell=(\ell_1,\ldots \ell_n)$ is a multi-index consisting of nonnegative
integers with order $|\ell|:=\ell_1+\ldots+\ell_n$, then we denote by $D^\ell$ the expression
$d_{e_1}^{\ell_1}\ldots d_{e_n}^{\ell_n}$, while for any non-negative integer $k$,
$D^k:=\{D^\ell\, |\, |\ell|=k\}$. For a point $p$ in, and a function $f$ on, the frame domain,
$D^k\!f|_p:=\{D^\ell\!f|_p\, |\, |\ell|=k\}$ can be regarded, after ordering the various
$d_{e_i}^m$'s in some fashion, as a point in $\mathbb{R}^{C^k_n}$, where $C^k_n$ is the number of
partitions of $k$ as an ordered sum of $n$ non-negative integers.

With these notations, by a {\em generalized PDE} (or {\em frame-dependent PDE}) of order $k$ on
a domain $U$ contained in the frame domain, we mean  an expression of the form
\be\lb{gen-pde} F(D^k\!f|_p, D^{k-1}\!f|_p,\ldots, D\!f|_p, f(p),p)=0 \end{equation}
where $f$ is an unknown standing for a function on $U$, whereas
\[
F:\mathbb{R}^{C^k_n}\times\mathbb{R}^{C_n^{k-1}}\times\ldots\mathbb{R}^{C^1_n}\times\mathbb{R}\times U\to\mathbb{R}
\]
is given. Slightly extending this definition, we also allow any occurrence of $f$ in \Ref{gen-pde}
to be replaced by $h\circ f$, where $h$ is a smooth real valued function defined on the image of $f$.


In our theorems the generalized PDE for the twist function is second order, and
in the case of Theorem~\ref{ke} it is equivalent to the requirement of constant
Gauss curvature of a related Riemmanian surface metric. It is also an
integral part of Theorem~\ref{KEc}. In Theorem~\ref{central}
this generalized PDE plays only an indirect role: central curvature zero turns out
to be a partial curvature condition effecting only an eigenspace of the
Ricci endomorphism on which the directional derivatives of $\iota$ are known,
so the generalized PDE, which involves frame fields of the {\em other} vector
fields, is not needed to determine centrality of the metric. It only appears in the
part of the theorem stating that the central metric is conformal to a CSC metric.



\section{Background on the construction of K\"ahler metrics in \cite{a-m}}
\label{sec:1}
Let $(M,g)$ be an oriented semi-Riemannian $4$-manifold with two vector fields
$\kk$, $\tT$ such that
\be\lb{basic0}
\begin{aligned}
&\text{$\kk$, $\tT$ are everywhere linearly independent;}\\
&\text{for $\VV:=\mathrm{span}(\kk,\tT)$, the distribution $\HH:=\VV^\perp$ is spacelike.}
\end{aligned}
\end{equation}
Here $\mathrm{span}(\kk,\tT)$ denotes the distribution spanned pointwise by
these vector fields, $\perp$ denotes the orthogonal complement, and the spacelike
requirement means that the restriction $g|_\HH$ of $g$ to $\HH$ is positive definite.
It also implies $TM=\HH\oplus\VV$ (see \cite[Lemmas 2.23, 2.22(2)]{on}) and that $\HH$ has rank two.
Assume that
\be\lb{integra}
[\kk,\Gamma(\HH)]\subset\Gamma(\HH),\qquad [\tT,\Gamma(\HH)]\subset\Gamma(\HH)
\end{equation}
and
\be\lb{shear}
J\n^o \kk=\n^o\tT
\end{equation}
where $J=J_{g,\kk,\tT}$ is the almost complex structure taking $\kk$ to $\tT$ and
making $g|_\HH$ hermitian and compatible with the orientation,
while $\n^o\!X$ is the shear operator of a vector field $X$:
the trace-free symmetric part of $\pi\circ\n X\big|_\HH:\HH\to\HH$,
for the orthogonal projection $\pi:TM\to\HH$. If $\n^o\!X=0$ we say $X$ is shear-free.

Assumptions \Ref{basic0} are necessary for $J$ to be well-defined.
In turn, by Theorem~$1$ of \cite{a-m}, assumptions \Ref{integra}
and \Ref{shear} guarantee that $J$ is integrable.

Further assumptions are needed for the existence of a class
of K\"ahler metrics on certain regions in $M$. These
are
\be\lb{kahler}
\begin{aligned}
&\text{$\tT=\ell\n\ta$, for $C^\infty$ functions $\ta$, $\ell$,}\\
&\n(g(\kk,\tT))\in\Gamma(\VV),\qquad \n(g(\kk,\kk))\in\Gamma(\VV).
\end{aligned}
\end{equation}
By Theorem~$3$ of \cite{a-m}, conditions \Ref{basic0}, \Ref{integra}, \Ref{shear} and \Ref{kahler}
guarantee the existence of K\"ahler metrics of the form
\be\lb{gen-ka}
\gK=-d(f(\ta) \kk^\flat)(J\cdot, \cdot)
\end{equation}
on any region satisfying
\be\lb{dom}
f\iota<0,\qquad f'\det(g|_\VV)/\ell-fd\kk^\flat(\kk,\tT)<0,
\end{equation}
where the notation is as follows. $f$ is a smooth real-valued function on a subset of $\mathbb{R}$,
$\kk^\flat$ denotes the $1$-form $g$-dual to $\kk$, and the prime denotes differentiation with
respect to $f$'s variable. Next, $\det(g|_\VV)=G:=g(\kk,\kk)g(\tT,\tT)-g(\kk,\tT)^2$.
Finally, $\imath:=g(\kk,[\xx,\yy])$ is the twist function
of $\kk$, for any (local) oriented orthonormal frame $\xx$, $\yy$ of $\HH$.
Here by an oriented frame, we mean one aligned with the orientation
induced on $\HH$ by the orientations of $M$ and $\VV$, the latter determined by
the ordered pair $\kk$, $\tT$. Within a given common domain for a class of such frames,
the twist function is independent of the choice of frame.

In the following sections we will often be strengthening assumptions \Ref{kahler}
and adding some additional Lie bracket related assumptions to \Ref{integra}.

We will occasionally need the expression for the shear relative to an oriented
orthonormal frame $\xx$, $\yy$ for $\HH$. At each point,
the matrix of the shear operator of a vector field $X\in\Gamma(\VV)$ is
\be\lb{shr0}
[\n^o\!X]_{\xx,\yy}=
\begin{bmatrix}
        - \sigma_1 & \sigma_2\\
        \sigma_2 & \sigma_1\\
      \end{bmatrix},
\end{equation}
where the entries are the \emph{shear coefficients}
\be
\begin{aligned}
\sigma_1\ &=\ 
\ \ \frac{1}{2}\Big[g([X,\xx],\xx)-g([X,\yy],\yy)\Big],\\
\sigma_2\ &=\ 
 -\frac{1}{2}\Big[g([X,\xx],\yy)+g([X,\yy],\xx)\Big]\cdot\label{eqn:shear2}
\end{aligned}
\end{equation}

We now give some definitions. The quadruple $(M,g,\kk,\tT)$ as above, satisfying
\Ref{basic0}, \Ref{integra}, \Ref{shear}, \Ref{kahler} is called {\em admissible},
and so is $g$ itself. The K\"ahler metric $\gK$ is said to be {\em induced} by the
admissible structure or metric.
We also use the terminology for vector fields used in \cite{a-m}:
$\kk$ is called {\em geodesic} if $\n_\kk\kk=0$,  {\em strictly pre-geodesic} if
$\n_\kk\kk=\gamma\kk$, $\gamma\ne 0$, and of course {\em Killing} if $\mathcal{L}_\kk g=0$,
where $\mathcal{L}$ is the Lie derivative.

In order to make the text more self-contained, we now list some facts whose proofs are in
\cite{a-m}. First, according to Remark 4.2 in \cite{a-m}, the expression in the second inequality in \Ref{dom}, which is just the value $-\gK(\kk,\kk)=-\gK(\tT,\tT)$, simplifies when $\kk$ is distinguished, so that $\gK(\kk,\kk)=\gK(\tT,\tT)=$
\begin{align}\lb{gkk}
&\text{$-f'G/\ell$ if $\kk$ is geodesic of constant length,}\nonumber\\
&\text{$-f'G/\ell + f\gamma g(\kk,\tT)$ if $\kk$ is null and strictly pre-geodesic,}\nonumber\\
&\text{$-f'G/\ell - fd_\tT(g(k,k))$ if $\kk$ is Killing.}
\end{align}

Next, the K\"ahler metric $\gK$ of the form \Ref{gen-ka}
has the following properties
\cite[Lemma~4.4]{a-m}:
\be\lb{gK}
\gK(\HH,\VV)=0, \qquad \gK(\kk,\tT)=0, \qquad \gK|_\HH=-f\iota g|_\HH.
\end{equation}

\begin{remark}\lb{shearK}
The shear operators of $\kk$, $\tT$ with respect to $\gK$ are equal to those with
respect to $g$ \cite[formula (25)]{a-m}. Note that expressions analogous to
\Ref{eqn:shear2} but with $\gK$ replacing $g$ will be studied later. They differ from
the actual shear coefficients of $\gK$ by a function multiple, because the last part of
\Ref{gK} shows that $\xx$, $\yy$ are only $\gK$-orthogonal, rather than $\gK$-orthonormal.
\end{remark}
\begin{remark}\lb{Nintegra}
Note that the proof of \cite[Theorem 1]{a-m} shows that to determine integrability
of $J$ it is enough to show that $N(\kk,\xx)=0$, where $N$ is the Nijenhaus tensor
and $\xx\in\Gamma(\HH)$.
\end{remark}
\begin{remark}\lb{twist-gr}
The twist function of a gradient vanishes, for example we have $g(\n\ta,[\xx,\yy])=
-g(\n_\xx\n\ta,\yy)+g(\n_\yy\n\ta,\xx)=0$ as the Hessian is symmetric.
\end{remark}

\section{Central and K\"ahler-Einstein metrics on admissible manifolds}
\lb{sec:2}
Aside from the stipulations on admissible manifolds described
in Section \ref{sec:1}, we will need in this section
additional Lie bracket requirements that will simplify our calculations.

Specifically, we require the existence, in a neighborhood of each point of $M$,
of an oriented orthonormal frame $\xx$, $\yy$  for $\HH$ satisfying
\be\lb{*}
\begin{aligned}
&[\kk,\xx]=\al\yy,\qquad &&[\tT,\xx]=\bet\yy,\quad\text {for constants $\al$, $\bet$},\\
&[\kk,\yy]\in\Gamma(\yy^\perp),\quad &&[\tT,\yy]\in\Gamma(\yy^\perp),\quad [\xx,\yy]\in\Gamma(\VV).
\end{aligned}
\end{equation}
Examples of admissible manifolds satisfying \Ref{*} will be given in Sections~\ref{GR} and \ref{Lie-g}.

Note that the second row of \Ref{*} represents relations that are somewhat weaker
than the requirements needed for $k$-dependency. This will turn out not
to cause any difficulty because the Koszul formula involves not Lie brackets per se,
but specific metric values on a Lie bracket and another vector field. Also,
in the course of the  following theorem's proof, exact Lie bracket relations will be deduced,
which conform to the $2$-dependency mentioned in subsection~\ref{k-dep}.

Also as in that subsection, in the following we denote by $d_X$ the derivative of a function
in the direction of a vector field $X$.
To state our first result, recall that the {\em central curvature} of a metric \cite{cen}
is the determinant of its Ricci endomorphism, and if it is constant for a K\"ahler
metric, the metric is also called {\em central}.
\begin{thm}\lb{central}
Let $(M,g)$ be an admissible $4$-manifold with $\kk$ either a geodesic vector field
or a Killing field. Assume additionally that $\kk$, $\tT$ commute, are
shear-free and satisfy conditions \Ref{*} for a local oriented
orthonormal frame $\xx$, $\yy$ of $\HH$ near each point.
Suppose $\kk$ is null, $\tT=\n\ta$ is gradient,
$g(\kk,\tT)=a$, $g(\tT,\tT)=b$ are constants
and $\n\iota\in\Gamma(\HH)$, where $\iota$ is
the twist function of $\kk$.
Then, wherever defined,  the K\"ahler metric
\[
\gK=-d(e^\ta \kk^\flat)(J\cdot, \cdot)
\]
is central, with vanishing central curvature.
Moreover, the conformally related metric $\tilde{g}=e^{-\ta}\!\gK$
is CSC precisely when $\iota$ satisfies the generalized PDE
\be\lb{pde}
(d_\xx d_\xx+d_\yy d_\yy)(\log|\iota|)=c\iota
\end{equation}
for some constant $c$. If $\iota$ is actually constant, $\tilde{g}$ is also
locally  isometric to a left-invariant metric on a  Lie group.
\end{thm}
In more detail, for the Ricci endomorphism of $\gK$ at each point $p$,
the nonzero tangent vectors in $\VV_p$ are eigenvectors for the eigenvalue zero.
Whereas for constant $\iota$, $\HH_p$ is contained in an eigenspace for an eigenvalue
which is a constant multiple $q$ of $e^{-\ta}$,
so that the scalar curvature of $\gK$ is $2qe^{-\ta}$.
Thus $\gK$ is either Ricci flat (in fact flat) or has semidefinite Ricci curvature
whose sign is that of $q$. Finally, if $\iota$ is constant, the
Levi-Civita connection of $\gK$ is also left invariant.
For general studies of invariant K\"ahler metrics on $4$-dimensional Lie groups
see \cite{ov}. A description of the possible Lie groups will be given in
subsection~\ref{Lie-g}.

The proof will be given in the next three subsections, followed by two subsections
describing examples. An analogous theorem will follow for K\"ahler-Einstein metrics.

\subsection{The connection}\lb{conn-ct}
As our metric is admissible, any induced K\"ahler metric $\gK$ of the form \Ref{gen-ka}
satisfies \Ref{gK} in its domain, while it follows from the fact that
$\kk$ is geodesic or Killing with $\kk$-null along with formulas \Ref{gkk} that
\be\lb{gk-formula}
\gK(\kk,\kk)=\gK(\tT,\tT)=-f'G/\ell=a^2f'.
\end{equation}
We remark that $a\ne 0$, as $g\big|_\VV$ is non-degenerate at each point for an admissible metric.

In the following we calculate in the frame $\kk$, $\tT$, $\xx$, $\yy$, with the last two vector fields defined locally and satisfying \Ref{*}.

We note that $d_\kk\ta=g(k,\n\ta)=a$ and similarly $d_\tT\ta=b$, while $d_\xx\ta=g(x,\n\ta)=0$ as $\n\ta\in\Gamma(\VV)$, and similarly $d_\yy\ta=0$. Therefore \Ref{gk-formula} implies
\[
d_\xx(\gK(\kk,\kk))=d_\yy(\gK(\tT,\tT))=0,\quad d_\kk(\gK(\kk,\kk))=a^3f'',\quad  d_\tT(\gK(\kk,\kk))=a^2bf''.
\]
By \Ref{gK}, our frame is $\gK$-orthogonal, so that
\be\lb{Kosz}
\nK_ab=\fr{\gK(\nK_ab,\xx)}{\gK(\xx,\xx)}\xx+\fr{\gK(\nK_ab,\yy)}{\gK(\yy,\yy)}\yy
+\fr{\gK(\nK_ab,\kk)}{\gK(\kk,\kk)}\kk+\fr{\gK(\nK_ab,\tT)}{\gK(\tT,\tT)}\tT.
\end{equation}
We apply this formula repeatedly in the following lemmas, using the Koszul formula
\[
\begin{aligned}
2\gK(\nK_ab,c)&=d_a(\gK(b,c))+d_b(\gK(a,c))-d_c(\gK(a,b))\\
&-\gK(a,[b,c])\ -\gK(b,[a,c])\,+\gK(c,[a,b])
\end{aligned}
\]
to compute the numerators in it. (Note that it is enough
to compute only three well-chosen covariant derivatives on frame fields
such as $\nK_\kk\tT$, $\nK_\xx\yy$, $\nK_\kk\xx$
to obtain all the other covariant derivatives on frame fields
via the requirements that $\nK$ be torsion-free and make $J_{g,\kk,\tT}$ parallel.)

For example, the metric values and their directional derivatives that
we have just computed, along with \Ref{integra}
and $[\kk,\tT]=0$ are used in the Koszul formula to give:
\begin{lemma}\lb{1}
\[
\begin{aligned}
\nK_\kk\kk&=-\nK_\tT\tT=\fr{f''}{2f'}(a\kk-b\tT),\qquad
\nK_\kk\tT=\nK_\tT\kk=\fr{f''}{2f'}(b\kk+a\tT).
\end{aligned}
\]
\end{lemma}

To obtain the next covariant derivative formulas we will be using
conditions \Ref{*}. In fact it is sufficient that the first
two of these conditions be weakened to $[\kk,\xx]\in\xx^\perp$,
$[\tT,\xx]\in\xx^\perp$.

Additionally, we will need
\be\lb{shear1}
\gK([\kk,\xx],\yy)+\gK([\kk,\yy],\xx)=0,
\end{equation}
along with the same relation with $\kk$ replaced by $\tT$.
These follow from the $g$-shear-freeness of $\kk$, $\tT$, the second
equation in \Ref{eqn:shear2} and Remark~\ref{shearK}.


Putting these relations and \Ref{*} in the Koszul formula
we arrive at
\be\lb{temp}
\begin{aligned}
\nK_\xx\kk&=\fr{d_\kk\gK(\xx,\xx)}{2\gK(\xx,\xx)}\xx
-\fr{\iota_{\scriptscriptstyle K}}{2\gK(\yy,\yy)}\yy,\quad
\nK_\yy\kk=\fr{\iota_{\scriptscriptstyle K}}{2\gK(\xx,\xx)}\xx
+\fr{d_\kk\gK(\yy,\yy)}{2\gK(\yy,\yy)}\yy,\\
\nK_\xx\tT&=\fr{d_\tT\gK(\xx,\xx)}{2\gK(\xx,\xx)}\xx
-\fr{\iota_{\scriptscriptstyle K}^\tT}{2\gK(\yy,\yy)}\yy,\quad
\nK_\yy\tT=\fr{\iota_{\scriptscriptstyle K}^\tT}{2\gK(\xx,\xx)}\xx
+\fr{d_\tT\gK(\yy,\yy)}{2\gK(\yy,\yy)}\yy,
\end{aligned}
\end{equation}
where we used the notations $\iota_{\scriptscriptstyle K}:=\gK(\kk,[\xx,\yy])$,
$\iota_{\scriptscriptstyle K}^\tT:=\gK(\tT,[\xx,\yy])$, even though these
quantities are again not quite the $\gK$-twists of $\kk$, $\tT$.

We now compute the numerators in the above expressions.
From the third relation in \Ref{gK}, the condition that
the twist of $\kk$ has a horizontal gradient, and the
frame directional derivatives of $\ta$ computed earlier, we see that
\[
\begin{aligned}
d_\kk\gK(\xx,\xx)&=d_\kk\gK(\yy,\yy)=-d_\kk(f\iota)=-\iota af',\\
d_\tT\gK(\xx,\xx)&=d_\tT\gK(\yy,\yy)=-d_\tT(f\iota)=-\iota bf'.
\end{aligned}
\]

We now compute $\iota_{\scriptscriptstyle K}$, $\iota_{\scriptscriptstyle K}^\tT$.
Using the last part of \Ref{*} we write $[\xx,\yy]=r\kk+s\tT$ for some coefficients $r$, $s$.
Then
\begin{align*}
\iota=g(\kk,[\xx,\yy])&=sa,\\
0=\iota^t=g(\tT,[\xx,\yy])&=ra+sb,\\
\iota_{\scriptscriptstyle K}=\gK(\kk,[\xx,\yy])&=ra^2f',\\
\iota_{\scriptscriptstyle K}^\tT=\gK(\tT,[\xx,\yy])&=sa^2f'.\\
\end{align*}
Here we have used the metric values of $g$ and $\gK$ on $\kk$, $\tT$
(see Theorem~\ref{central}'s statement and \Ref{gk-formula}). The fact that $\iota^t=0$
follows from the fact that $\tT$ is gradient and Remark~\ref{twist-gr}.

Substituting the first of these four equations in the last, and the second in the third
gives $\iota_{\scriptscriptstyle K}=-sbaf'$ and
$\iota_{\scriptscriptstyle K}^\tT=\iota af'$. Replacing again $sa$ by $\iota$ in the first of these
finally yields
\[
\iota_{\scriptscriptstyle K}=-\iota bf'\qquad    \iota_{\scriptscriptstyle K}^\tT=\iota af'.
\]
Substituting the above expressions in the numerators of \Ref{temp}, and the denominators using
the third relation in \Ref{gK}, we finally get
\begin{lemma}\lb{2}
\[
\begin{aligned}
\nK_\xx\kk&=\fr{f'}{2f}(a\xx-b\yy),
\qquad
\nK_\yy\kk=\fr{f'}{2f}(b\xx+a\yy),\\
\nK_\xx\tT&=\fr{f'}{2f}(b\xx+a\yy),
\qquad
\nK_\yy\tT=\fr {f'}{2f}(-a\xx+b\yy).
\end{aligned}
\]
\end{lemma}
When applying the above formulas for
$\iota_{\scriptscriptstyle K}$,  $\iota_{\scriptscriptstyle K}^\tT$
together with \Ref{*}, conditions \Ref{shear1}
(for $\tT$ as well) and the fact that $\ta$ has a vertical gradient,
one arrives similarly at
\begin{lemma}\lb{3}
\[
\begin{aligned}
\nK_\xx\yy&=\fr 1{2\iota}(d_\yy\iota\,\xx+d_\xx\iota\,\yy)+\fr{\iota}{2a^2}\left(- b\kk+ a\tT\right),\\
\nK_\yy\xx&=\fr 1{2\iota}(d_\yy\iota\,\xx+d_\xx\iota\,\yy)-\fr{\iota}{2a^2}\left(- b\kk+a\tT\right),\\
\nK_\xx\xx&=\fr 1{2\iota}(d_\xx\iota\,\xx-d_\yy\iota\,\yy)+\fr{\iota}{2a^2}\left(a \kk+b\tT\right),\\
\nK_\yy\yy&=\fr 1{2\iota}(-d_\xx\iota\,\xx+d_\yy\iota\,\yy)+\fr{\iota}{2a^2}\left(a \kk+b\tT\right).
\end{aligned}
\]
\end{lemma}

Finally, we note
\begin{lemma}\lb{4}
\[
\begin{aligned}
\nK_\kk\xx&=\nK_\xx\kk+\al\yy,\qquad\nK_\tT\xx&&=\nK_\xx\tT+\bet\yy,\\
\nK_\kk\yy&=\nK_\yy\kk-\al\xx,\qquad\nK_\tT\yy&&=\nK_\yy\tT-\bet\xx.
\end{aligned}
\]
\end{lemma}
The first line of this lemma follows from the fact that $\nK$ is torsion-free
and from the first line of \Ref{*}. The lemma's second line follows similarly,
because \Ref{*}, \Ref{integra} and the shear-freeness of $\kk$ and $\tT$
(expressed via its shear coefficients as in the second line of \Ref{eqn:shear2}) imply
\be\lb{brack-part}
[\kk,\yy]=-\al\xx,\qquad [\tT,\yy]=-\bet\xx.
\end{equation}

We combine the information in the last four lemmas as follows. Consider $(TM, J)$
as a complex bundle and set $w_1:=\kk-i\tT$, $w_2:=\xx-i\yy$. We use the same notation $\nK$
for the complexified connection obtained by extending linearly the Levi-Civita connection
$\nK$ of $\gK$, so that it differentiates complex vectors fields such as the $w_i$'s
along {\em real directions}, for example in the directions of our standard (non-complex)
frame vector fields. Then, using the summation convention, we write
\be\lb{compl-conn}
\nK w_i=\Gamma_i^j\otimes w_j,
\end{equation}
where the $\Gamma_i^j$ are complex-valued $1$-forms, whose expression we can compute
by applying the above 4 lemmas. Specifically, $\Gamma_i^j=\Gamma_i^j(\xx)\hat{\xx}+\Gamma_i^j(\yy)\hat{\yy}+
\Gamma_i^j(\kk)\hat{\kk}+\Gamma_i^j(\tT)\hat{\tT}$, where the hatted quantities constitute
the non-metric dual coframe to our frame, so that, for example, $\hat{\xx}$ is $1$ on $\xx$ and zero
on the other frame fields. The coefficients in this expansion are calculated by
substituting in \Ref{compl-conn} the frame vector fields and then using the
expressions in Lemmas \ref{1}, \ref{2}, \ref{3} and \ref{4}.
The final result of these calculations is given by
\be\lb{gammas}
\begin{aligned}
\Gamma_1^1&=\fr{f''}{2f'}(a-ib)\hat{\kk}+\fr{f''}{2f'}(b+ia)\hat{\tT},\\
\Gamma_1^2&=\fr{f'}{2f}(a-ib)\hat{\xx}+\fr{f'}{2f}(b+ia)\hat{\yy},\\
\Gamma_2^1&=\fr{\iota}{2a^2}(a+ib)\hat{\xx}+\fr{\iota}{2a^2}(b-ia)\hat{\yy},\\
\Gamma_2^2&=\fr 1{2\iota}(d_\xx\iota-id_\yy\iota)\hat{\xx}+
\fr 1{2\iota}(d_\yy\iota+id_\xx\iota)\hat{\yy}\\
&+\left[\fr{f'}{2f}(a-ib)+\alpha i\right]\!\hat{\kk}+
\left[\fr{f'}{2f}(b+ia)+\beta i\right]\!\hat{\tT}.
\end{aligned}
\end{equation}

\subsection{The Ricci curvature}
According to Lemma $4.2$ of \cite{dm1}, the Ricci form of $\gK$ is given
by
\be\lb{rho}
\rho_{\scriptscriptstyle K}=i(d\Gamma_1^1+d\Gamma_2^2).
\end{equation}
We now wish to evaluate it on our frame vector fields.
From now on we fix $f(\ta)=e^\ta$.

Due to the formula
\be\lb{1-form}
d\xi(u,v)=d_u(\xi(v))-d_v(\xi(u))-\xi([u,v]),
\end{equation}
applied to the frame vector fields and the coframe $1$-forms,
we see that $d\hat{\kk}$, $d\hat{\tT}$, $d\hat{\xx}$, $d\hat{\yy}$
vanish on vector fields in $\VV$, since the coframe $1$-forms have constant
values on the frame fields while $\kk$, $\tT$ commute.
Given that $f(\ta)=e^\ta$, many of the coefficients of the connection
1-forms are constant. Those that may not be are the coefficients of $\hat{\xx}$, $\hat{\yy}$
in $\Gamma_2^2$. But as $\hat{\xx}$, $\hat{\yy}$ actually vanish on $\kk$, $\tT$,
it follows that $d\Gamma_1^1$, $d\Gamma_2^2$ and therefore $\rho_{\scriptscriptstyle K}$
all vanish on $\kk$, $\tT$.

We now show the the Ricci form vanishes for a pair of frame fields, one in $\VV$
and the other in $\HH$. The argument is entirely analogous to the above
except for the terms of $\Gamma_2^2$ involving $\hat{\xx}$ and $\hat{\yy}$.
To see what happens in those, consider for example the case where the frame fields
are $\kk$ and $\xx$. The first relation in \Ref{*}, the values of the coframe
on the frame and the fact that $\n\iota$ lies in $\HH$, so $d_\kk\iota=0$, mean
that there are only two non-vanishing terms in $d\Gamma_2^2$. The first is
$d[(1/(2\iota))(d_\xx\iota-id_\yy\iota)](\kk)\hat{\xx}(\xx)=
(1/(2\iota))(d_\kk d_\xx\iota-id_\kk d_\yy\iota)=
(1/(2\iota))(\al d_\yy\iota+i\al d_\xx\iota)$, where the last
equality uses the first line of \Ref{*} and the vanishing of $d_\kk\iota$.
The second non-vanishing term is
$(1/(2\iota))(d_\yy\iota+id_\xx\iota)d\hat{\yy}(\kk,\xx)=
(1/(2\iota))(d_\yy\iota+id_\xx\iota)(-\hat{\yy}([\kk,\xx])
=(1/(2\iota))(d_\yy\iota+id_\xx\iota)(-\al)$, using \Ref{*}. We thus see
that these two terms cancel each other, proving the claim for $\kk$ and $\xx$.
For other pairs, such that one vector field is in $\VV$ and the other is in $\HH$,
the proof is similar.

Thus the Ricci tensor of $\gK$ has a zero eigenvalue, with eigenvectors which
include the nonzero tangent vectors of $\VV$ at a given point. Hence $\gK$ is central.

Note that by Lemma \ref{3}
\be\lb{***}
\text{$[\xx,\yy]=\fr\iota{a^2}(-b\kk+a\tT)$,}
\end{equation}
because $[\xx,\yy]=\nK_\xx\yy-\nK_\yy\xx$.
This and \Ref{1-form} give $d\hat{\kk}(\xx,\yy)=-\hat{\kk}([\xx,\yy])
=-\hat{\kk}((\iota/a^2)(-b\kk+a\tT))=b\iota/a^2$,
and similarly $d\hat{\tT}(\xx,\yy)=-\iota/a$.
We thus compute via \Ref{rho}
\[
\begin{aligned}
\rho_{\scriptscriptstyle K}(\xx,\yy)&=\fr12 i(b\iota/a^2)(a-ib)+\fr12 i(-a\iota/a^2)(b+ia)\\
&+\fr12 i(b\iota/a^2)(a-ib)+\fr12 i(-a\iota/a^2)(b+ia)
+i^2\alpha b\iota/a^2-i^2\beta a\iota/a^2\\
&-id\left[\fr 1{2\iota}(d_\xx\iota-id_\yy\iota)\right](\yy)
+id\left[\fr 1{2\iota}(d_\yy\iota+id_\xx\iota)\right](\xx)\\
&=\fr\iota{a^2}(a^2+b^2)+\iota\fr{-b\alpha+a\beta}{a^2}\\
&-\fr i2\left(-\fr {d_\yy\iota}{\iota^2}(d_\xx\iota-id_\yy\iota)+
\fr1\iota(d_\yy d_\xx\iota-id_\yy d_\yy\iota)\right)\\
&+\fr i2\left(-\fr {d_\xx\iota}{\iota^2}(d_\yy\iota+id_\xx\iota)+
\fr 1\iota(d_\xx d_\yy\iota+id_\xx d_\xx\iota)\right)\\
&=\iota\fr{a^2+b^2-b\alpha+a\beta}{a^2}\\
&+\fr 1{2\iota^2}((d_\xx\iota)^2+(d_\yy\iota)^2)
-\fr 1{2\iota}(d_\xx d_\xx\iota+d_\yy d_\yy\iota)\\
&=\iota\fr{a^2+b^2-b\alpha+a\beta}{a^2}-
\fr12(d_\xx d_\xx+d_\yy d_\yy)(\log|\iota|),
\end{aligned}
\]
where in the next to last step we used the fact that $d_\yy d_\xx\iota-d_\xx d_\yy\iota=0$,
which follows from the above formula for $[\xx,\yy]$ and the assumption that
$\n\iota$ has a horizontal gradient. This calculation will be needed later in the proof.

Assume now that $\iota$ is constant. Then the above calculation shows that
the Ricci curvature of $\gK$ has a constant value on the pair $\xx$, $\xx$ or on $\yy$, $\yy$.
If this constant Ricci value is denoted $\ell$, by the third relation in \Ref{gK}
$-\ell/(\iota e^\ta)$ is an eigenvalue of the Ricci tensor with eigenspace containing $\HH$.
The constant $q$ appearing in the paragraph after the statement of Theorem~\ref{central} is thus $q=-\ell/\iota=-(a^2+b^2-b\alpha+a\beta)/a^2$.

Finally, note that if $\iota$ is constant, the covariant derivatives have constant coefficients in the frame $\kk,\tT,\xx,\yy$. By the fact that $\kk$, $\tT$ commute as well as the Lie bracket expressions
for our frame given in \Ref{*}, \Ref{brack-part} and \Ref{***}, this frame generates a Lie algebra. Hence the Levi-Civita connection is locally a {\em left-invariant} torsion-free connection for any Lie group
whose tangent space at the identity is this Lie algebra.

\begin{remark}
One can check that the $\gK$-curvature tensor values on our
frame fields all vanish except on $\xx$, $\yy$, so that $\gK$ is
Ricci-flat if and only if it is flat.
\end{remark}

\subsection{The conformally related metric}
We now turn to the conformally related metric $\tilde{g}=e^{-\ta}\!g_K$.
From general conformal change formulas, its scalar curvature is
\be\lb{conf-scal}
\tilde{s}=s_{\!\scriptscriptstyle K}u^2+
6u\Delta^{\!\scriptscriptstyle K}\! u-12\gK\!(\nK\! u,\nK\! u),\qquad
\text{for }u=e^{\ta/2},
\end{equation}
where $s_{\scriptscriptstyle K}$ is the scalar curvature of $\gK$,
$\Delta^{\!\scriptscriptstyle K}$ the $\gK$-Laplacian
and $\nK\!u$ the $\gK$-gradient of $u$.

We now demonstrate that all three summands in this formula are constant,
beginning with the case where {\em $\iota$ is constant}.

First, $s_{\scriptscriptstyle K}u^2=2qe^{-\ta}e^\ta=2q$ is constant.
Next, $\nK\! u=\gK^{-1}(du, \cdot)=(e^{\ta/2}/2)\gK^{-1}(d\ta,\cdot)
=(e^{\ta/2}/2)\gK^{-1}(a\hat{\kk}+b\hat{\tT},\cdot)=
(e^{\ta/2}/(2\gK(\kk,\kk))(a\kk+b\tT)$,
so that by our formulas for $\gK(\kk,\kk)=\gK(\tT,\tT)$ we see
that  $\gK(\nK\!u,\nK\! u)=e^\ta(a^2+b^2)a^2e^\ta/(4a^4e^{2\ta})=
(a^2+b^2)/(4a^2)$ is also constant.

Finally, let $\{e_i\}$ be the $\gK$-orthonormal frame obtained from our
standard frame. We wish to apply the Laplacian formula
\[
\Delta^{\scriptscriptstyle K}\!\ta=\sum_i(d_{e_i}^2\ta-d\ta(\nK_{e_i}e_i)).
\]
The first term of each summand in this formula contributes only for $e_1=\kk/|\kk|$ and $e_2=\tT/|\tT|$
(with $|\cdot|$ denoting the $\gK$-norm), as $d_\xx\ta=d_\yy\ta=0$. The first of these
contributions is
$d_{\kk/|\kk|}d_{\kk/|\kk|}\ta=a|\kk|^{-1}(-|\kk|^{-2}d_\kk|\kk|)=-a|\kk|^{-3}(a^3e^{\ta}/(2|\kk|)
=-e^{-\ta}/2$. The second gives similarly $d_{\tT/|\tT|}d_{\tT/|\tT|}\ta=-(b^2/(2a^2))e^{-\ta}$.

Next, we compute the second part of each summand in the Laplacian formula, using Lemma \ref{1}
\[
\begin{aligned}
\nK_{\kk/|\kk|}(\kk/|\kk|)&=|\kk|^{-1}[d_\kk(|\kk|^{-1})\kk+|\kk|^{-1}\nK_\kk\kk]\\[3pt]
&=|\kk|^{-1}[-|\kk|^{-2}d_\kk|\kk|]\kk+(|\kk|^{-2}/2)(a\kk-b\tT)\\[3pt]
&=-|\kk|^{-3}a^3e^\ta(2|\kk|)^{-1}\kk+(|\kk|^{-2}/2)(a\kk-b\tT)\\[3pt]
&=-(e^{-\ta}/(2a))\kk+(a^{-2}e^{-\ta}/2)(a\kk-b\tT).
\end{aligned}
\]
Similarly, one computes via Lemmas \ref{1} and \ref{3}
\[
\begin{aligned}
\nK_{\tT/|\tT|}(\tT/|\tT|)&=-(be^{-\ta}/(2a^2))\tT-(a^{-2}e^{-\ta}/2)(a\kk-b\tT),\\[3pt]
\nK_{\xx/|\xx|}(\xx/|\xx|)&=\nK_{\yy/|\yy|}(\yy/|\yy|)=-(e^{-\ta}/2)(\fr 1a\kk+\fr b{a^2}\tT).
\end{aligned}
\]
Gathering all these partial calculations for the summands thus yields
\[
\begin{aligned}
\Delta^{\scriptscriptstyle K}\ta&=e^{-\ta}(-1/2-b^2/(2a^2))\\
&-e^{-\ta}[-1/2+1/2-b^2/(2a^2)-b^2/(2a^2)-1/2+b^2/(2a^2)-1-b^2/a^2]\\
&=-e^{-\ta}[-1-b^2/a^2]:=-e^{-\ta}v
\end{aligned}
\]
Now $\nK du=(e^{\ta/2}/2)\nK d\ta+(e^{\ta/2}/4)d\ta\otimes d\ta$,
so
\[
\begin{aligned}
\Delta^{\scriptscriptstyle K}u&=(e^{\ta/2}/2)\Delta^{\scriptscriptstyle K}\ta+
(e^{\ta/2}/4)\sum_i(d\ta(e_i))^2\\
&=(e^{\ta/2}/2)(-e^{-\ta}v)+(e^{\ta/2}/4)[a^2/(a^2e^{\ta})+b^2/(a^2e^{\ta})]\\
&:=e^{-\ta/2}p,
\end{aligned}
\]
where $p$ is a constant. Thus $6u\Delta^{\scriptscriptstyle K}u=6e^{\ta/2}e^{-\ta/2}p=6p$
is constant, hence so is the scalar curvature of $\tilde{g}$.
Putting this all together, the value of this scalar curvature is
\[
\begin{aligned}
\tilde{s}&=2q-12(a^2+b^2)/(4a^2)+6p\\
&=-2(a^2+b^2-b\alpha+a\beta)/a^2-3(a^2+b^2)/a^2\\
&+6(a^2+b^2)/(2a^2)+6(a^2+b^2)/(4a^2)\\
&=-\fr{a^2+b^2}{2a^2}+2\fr{b\alpha-a\beta}{a^2}.
\end{aligned}
\]

Recall that when $\iota$ is constant, our frame gives rise to a Lie algebra.
The formulas for $\gK$ on the frame, given partly in \Ref{gk-formula}, and
partly deduced from \Ref{gK} (for $f(\ta)=e^\ta$), show that $\tilde{g}=e^{-\ta}\gK$
has constant values on our frame when $\iota$ is constant. This implies $\tilde{g}$
is locally isometric to a left invariant metric with respect to a corresponding Lie group.

Finally we turn to the computation of $\tilde{s}$ when {\em $\iota$ is nonconstant}.
The third term on the right hand side of \Ref{conf-scal} yields the same constant
as before, by the same argument.
So does the second term there, since the additional terms in the computation of $\nK_{\xx/|\xx|}(\xx/|\xx|)$ and $\nK_{\yy/|\yy|}(\yy/|\yy|)$ are function combinations of $\xx$ and $\yy$, which vanish upon evaluation by $d\ta$.

This leaves the first term in \Ref{conf-scal}.
Now $s_{\scriptscriptstyle K}=2\rho_{\scriptscriptstyle K}(\xx,\yy)/\gK(\xx,\xx)
=2\rho_{\scriptscriptstyle K}(\xx,\yy)/(-e^\ta\iota)$. From the formula
for $\rho_{\scriptscriptstyle K}(\xx,\yy)$, given in the previous subsection,
we see that $s_{\scriptscriptstyle K}u^2$, and hence $\tilde{s}$,
will be constant exactly when the PDE \Ref{pde} holds for $\iota$ and some constant $c$.

This completes the proof of Theorem~\ref{central}.

\subsection{Examples related to General Relativity}\lb{GR}
We give two examples of this theorem with constant twist function $\iota$,
and another family with either constant or non-constant twist. Many features of these
are described in Sections~$8$ and $9$ of \cite{a-m}. The first example gives a trivial
K\"ahler metric but serves to check our formulas.
\begin{itemize}
\item{${\mathbf S^3\times\mathbb{R}}$}
Let $g$ be a Lorentzian product of the canonical Riemannian metrics on the $3$-sphere
and the real line, the latter equipped with coordinate $\ta$. The 3-sphere possesses an orthonormal frame $\bar{\kk}$, $\xx$, $\yy$ whose Lie brackets are given by cyclic permutations
of the relation $[\bar{\kk},\xx]=-2\yy$. On our $4$-manifold we choose our frame so that
$\xx$, $\yy$ are the lifts of the corresponding vector fields on the $3$-sphere, $\tT=\partial_\ta$ and $\kk=\bar{\kk}-\tT$. With these choices $g$ is admissible
with $\kk$ geodesic, $\kk$, $\tT$ commute and are shear-free, conditions \Ref{*} hold with
$\alpha=-2$ and $\beta=0$, $\kk$ is null, $\tT$ is the gradient of $\ta$, $a=-b=1$ and the twist
function $\iota$ is $-2$. By \cite[Theorem~5]{a-m}, the associated K\"ahler metric of
Theorem~\ref{central} is defined on the entire $4$-manifold, and the latter theorem shows it is
central and conformal to a (locally left invariant) constant scalar curvature metric. However the
constant defined when discussing $\rho_{\scriptscriptstyle{K}}(\xx,\yy)$ is
$q=-(a^2+b^2-b\al+a\bet)/a^2=0$, so that our metric is Ricci flat. Computation
of the sectional curvature from our covariant derivative formulas show that the
$\gK$ is in fact flat.

\item{{\bf Gravitational plane wave}}
This is a metric on $\mathbb{R}^4$ given by $g = -(x^2+y^2)du \otimes du + du \otimes dv
+ dv \otimes du + dx \otimes dx + dy \otimes dy$. The frame is $\kk=-\partial_u -y\partial_x + x\partial_y$, $\tT=\partial_v$, $\xx=-y\partial_v+\partial_x$, $\yy=x\partial_v+\partial_y$.
Then $g$ is admissible with $\kk$ geodesic (see \cite[Proposition 9.1]{a-m} and its proof),
$\kk$, $\tT$ commute and are shear-free, conditions \Ref{*} hold with
$\al=-1$, $\bet=0$, $\kk$ is null, $\tT=\n u$, $a=-1$, $b=0$, and $\iota=-2$.
The central K\"ahler metric $\gK$ in this case is defined on all of $\mathbb{R}^4$ and satisfies
$\mathrm{Ric}_{\scriptscriptstyle K}(\xx,\xx)=\iota=-2$.
Since $\iota$ is constant the metric $\tilde{g}=e^{-u}\gK$ is both CSC and locally left-invariant.

\item{{\bf Truncated\,pp-wave}}
This family of metrics $g$ is defined as the Lorentzian product of the standard metric
on the real line, equipped with coordinate $\ta$, and a metric $\bar{g}$ on $\mathbb{R}^3$
defined as follows. Let $h = H(u,x,y)du \otimes du + du \otimes dv + dv \otimes du + dx \otimes dx + dy \otimes dy$ be a pp-wave metric on $\mathbb{R}^4$, where $H$ is a smooth function.
Let $\xx = k(u,x,y)\partial_v + \partial_x$, $\yy = h(u,x,y)\partial_v + \partial_y$,
$\bar{\kk}=\partial_v$ for two smooth functions $k$, $h$ to be determined below. Restrict these
vector fields to a fixed $u=u_0$ hyperplane $S_{u_0}$ and define $\bar{g}$ by giving an
orthonormal coframe for it. Specifically, require $\bar{\xx}^{\bar\flat}$ (or
$\bar{\yy}^{\bar\flat}$) to be the restriction of $h(\xx,\cdot)$ (or $h(\yy,\cdot)$) to $S_{u_0}$, and $\bar{\kk}^{\bar\flat}$ a similar restriction of $-h(\zz,\cdot)$, where $\zz=\fr12(H + k^2 + h^2)\partial_v -\partial_u + k\partial_x + h\partial_y$. Our frame for the product metric then consists of
$\tT=-\partial_\ta$, $\kk=\bar{\kk}-\tT$, $\xx$, $\yy$.

Such a metric $g=g_{k,h}$ (which turns out not to depend on $H$), is admissible with $\kk$ a geodesic
vector field (\cite[Proposition 8.3 and Theorem 5]{a-m}).
The vector fields $\kk$, $\tT$ commute and are shear-free, \Ref{*} holds with $\al=\bet=0$,
$\kk$ is null, $\tT=\n\ta$, $a=1$ and $b=-1$. The central K\"ahler
metric is defined on all of $\mathbb{R}^4$ if the twist function of $\kk$,
which is $\iota=h_x-k_y$, is nowhere vanishing.

Now for the central K\"ahler metric $\gK$ we have
$\mathrm{Ric}_{\scriptscriptstyle K}(\xx,\xx)=2\iota-(1/2)(d_\xx d_\xx+d_\yy d_\yy)(\log(|\iota|)$.
By an appropriate choice of $k$, $h$, one could choose $\iota$ to be a (negative) constant as before, but there are other possible choices that will guarantee $\tilde{g}=e^{-\ta}\gK$ is CSC.

For example, choose $|\iota|=e^{p(x,y)}$ with $p$ a harmonic function in the $xy$-plane.
Since $d_\xx d_\xx+d_\yy d_\yy$ acts as the classical plane Laplacian on functions of $x$ and $y$, we see that $\iota$ will be a solution of \Ref{pde} with $c=0$, so the scalar curvature
$\tilde{s}$ of $\tilde{g}$ will have the same value $-1$ that it would have with the choice of
constant $\iota$.

Another possibility is to choose $\iota$ to be of the form $-\mathrm{sech}^2(z)$, where $z$ is affine
in $x$ and $y$ (and $c\ne 0$).
For some choice of the coefficients of $z$ one can guarantee, for example, $\tilde{s}=0$.


\end{itemize}

\subsection{Examples on Lie groups}\lb{Lie-g}

Consider the Lie bracket relations for our frame:
\begin{align}\lb{Lie-alg}
[\kk,\xx]&=\al\yy,  &&[\kk,\yy]=-\al\xx,\nonumber\\
[\tT,\xx]&=\bet\yy, &&[\tT,\,\yy]=-\bet\xx,\nonumber\\
[\kk,\tT]&=0,       &&[\xx,\yy]=(\iota/a^2)(-b\kk+a\tT).
\end{align}
Assume throughout this subsection that $\iota<0$ is constant.
We first address the question: which
Lie algebras $\mathfrak{g}$ are realized by these relations, up
to isomorphism?

Consider first cases when some of the constants in \Ref{Lie-alg}
vanish. Recall that $a\ne 0$ and $\iota\ne 0$. If $\al=\bet=0$,
then regardless of the value of $b$, the Lie algebra is $\mathfrak{nil}_3\times\mathbb{R}$,
where $\mathfrak{nil}_3$ is the Heisenberg Lie algebra.
This Lie algebra is realized by the vector fields in the truncated
pp-wave model (for constant $\iota$).
Next, if $\al=0$, $b=0$ and $\bet\ne 0$, then
$\mathfrak{g}=\mathfrak{su}(2)\times\mathbb{R}$
or $\mathfrak{g}=\mathfrak{sl}(2)\times\mathbb{R}$,
depending on the signs of $a$ and $\bet$. This will
also hold if $b\ne 0$, and either $\al=0$, $\bet\ne 0$ or
$\al\ne 0$, $\bet=0$ hold. On the other
hand, if $\al\ne 0$, $\bet=0$ and $b=0$ the Lie algebra
is unimodular, in fact step $3$ nilpotent,
classified in \cite{mc} as U3I2 (see case A8 in \cite{i-j-l}),
and does not correspond to any homogeneous compact geometry.
This Lie algebra is realized by the vector fields in our plane wave model.

Finally, in the generic case where $\al$, $\bet$, $a$, $b$
are all nonzero, again $\mathfrak{g}$ is either
$\mathfrak{su}(2)\times\mathbb{R}$
or U3I2, and the latter happens exactly when $-b\al+a\bet=0$. This
latter claim can be seen by the following argument communicated to us by
A. Derdzinski. Rescale $\kk$ to $\tilde{\kk}=\kk/\al$,
and $\tT$ to $\tilde{\tT}=\tT/\bet$, which makes
$[\xx,\yy]=(\iota/a^2)(-b\al\tilde{\kk}+a\bet\tilde{\tT})$.
Then if $-b\al+a\bet\ne 0$, rescale $\xx$, $\yy$ by the same
factor $p$ such that $p^3(\iota/a^2)(-b\al+a\bet)=:-p^3\iota B=1$
(the quantity $B$ will play a role later in Theorem~\ref{KEc}). Then choose
$\mathbf{p} = p^2(\iota/a^2)(-b\al\tilde{\kk} + a\bet\tilde{\tT})$ and
$\mathbf{q} = \tilde{\kk} - \tilde{\tT}$, so that
$\mathbf{p}$, $\mathbf{q}$, and the new $\xx$, $\yy$
satisfy the usual relations of $\mathbf{su}(2)\times\mathbb{R}$.
In the remaining case where $-b\al+a\bet= 0$, use the basis consisting of the original
$\xx$, $\yy$, then $\tilde{\kk}$ and $\mathbf{q} = (\iota/a^2)(-b\al\tilde{\kk} + a\bet\tilde{\tT})
= (\iota/a^2)(-b\al)(\tilde{\kk} - \tilde{\tT})$, which gives the
Lie brackets relations of U3I2, the non-zero ones being
$[\tilde\kk,\xx]=\yy$, $[\tilde\kk,\yy]=-\xx$ and $[\xx,\yy]=q$.

We consider now the fact that the constants $\alpha$, $\beta$ $a$, $b$ and $\iota$
determine not only the Lie algebra, but also, along with $f(\tau)$, our K\"ahler metric,
in a manner that is independent of the particular Lorentzian metric chosen (satisfying
Theorem~\ref{central}'s assumptions), so long as its value on $\kk$, $\tT$ is the given
constant $a$, its value on $\tT$, $\tT$ is the given $b$, and its twist is the given
$\iota$. This is because the symplectic form, when evaluated on frame fields, depends
on the Lorentzian metric only through quantities such as $\kk^\flat(v)$ or $\kk^\flat([v,w])$,
where $v$, $w$ are frame fields, but the theorem's assumptions fix these frame
values to be
$\kk^\flat(\kk)=\kk^\flat(\xx)=\kk^\flat(\yy)=0$, $\kk^\flat(\tT)=a$.
This observation will now allow us to carry out our construction globally on the
connected and simply connected Lie groups corresponding to the above four Lie algebras.

Let $\mathfrak{g}$ be a Lie algebra with a basis $\kk$, $\tT$, $\xx$, $\yy$,
satisfying relations \Ref{Lie-alg} for some constants
$\iota<0$, $\al$, $\bet$, $b$ and $a\ne 0$. Let these basis notations
also stand for the corresponding left-invariant vector
fields on the associated connected and simply-connected
Lie group $\mathcal{G}$. Define a semi-Riemannian
and left-invariant metric $g$ on $\mathcal{G}$
by $g(\xx,\xx)=g(\yy,\yy)=1$, $g(\xx,\yy)=g(\kk,\kk)=0$,
$g(\kk,\tT)=a$, $g(\tT,\tT)=b$ and $g(\VV,\HH)=0$,
where $\VV:=\mathrm{span}(\kk,\tT)$ and $\HH:=\mathrm{span}(\xx,\yy)$.
We now check that $g$ satisfies all conditions of
Theorem~\ref{central}. Clearly $g|_\HH$ is positive definite.
Relations \Ref{integra} follow from \Ref{Lie-alg}.
The second line of \Ref{kahler} is immediate as $\kk$ is
null and $a$ is constant. Similarly $\n\iota=0\in\Gamma(\HH)$.
One checks that $\kk$ and $\tT$ are shear-free by showing
via \Ref{Lie-alg} that their shear-coefficients
vanish, using formulas \Ref{eqn:shear2}.

It remains to show that $\kk$ is geodesic
and $\tT$ a gradient (and thus $g$ is admissible).
For this we use the following formula valid
for the Levi-Civita connection of a left-invariant
semi-Riemannian metric on a Lie group acting on
left-invariant vector fields (see \cite[Prop. 3.18]{ch-eb}):
\be\lb{Lie-conn}
\n_ab=\big([a,b]-\mathrm{ad}_a^*(b)-\mathrm{ad}_b^*(a)\big)/2,
\end{equation}
where $\mathrm{ad}_a^*$ is the adjoint of the $\mathrm{ad}$
operation. Hence $g(\n_\kk\kk,\cdot)=-g(\kk,[\kk,\cdot])$, and applying
this formula to the frame vector fields shows $\n_\kk\kk=0$. Next,
using
\Ref{Lie-conn} in the exterior derivative formula
$d\tT^\flat(a,b)=g(\n_a\tT,b)-g(\n_b\tT,a)$ of the $1$-form
$g$-dual to $\tT$ yields
$d\tT^\flat(a,b)=-g(\tT,[a,b])$. Applying this to any
pair of frame fields shows that $\tT^\flat$ is closed,
hence exact, as $\mathcal{G}$ is simply-connected. Thus
$\tT$ is a gradient, i.e. $\tT=\n\ta$ for some function $\ta$
on $\mathcal{G}$.

By Theorem~\ref{central} there is thus an induced central
K\"ahler metric $\gK$, whose domain is computed from \Ref{dom}
to be the entire group $\mathcal{G}$.
Hence (by the K\"ahler metric frame values in \Ref{gk-formula} and \Ref{gK},
for $f(\ta)=e^\ta$)
\[
\tilde{g}=e^{-\ta}\gK=\iota(\hat{\xx}^2+\hat{\yy}^2)+a^2(\hat{\kk}^2+\hat{\tT}^2),
\]
with the hatted quantities the usual non-metric duals to our frame,
is a (global) left-invariant metric (of course of CSC) on $\mathcal{G}$, hence complete.
Therefore $\gK$ is conformal to a complete CSC metric, but the unbounded conformal factor
$e^\ta$ prevents it from being complete (cf.~\cite{dps1,dps2}).

\subsection{K\"ahler-Einstein metrics}

We now obtain an analogous theorem for K\"ahler-Einstein metrics.

\begin{thm}\lb{KEc}
Let $(M,g)$ and $\kk$, $\tT$ satisfy all conditions of Theorem~\ref{central},
with $\iota$ satisfying the generalized PDE \Ref{pde} for a constant $c$.
Then, for any solution $f(\ta)$ of the ODE
\be\lb{KEceq} A\frac{(f^2)''(\ta)}{(f^2)'(\ta)}-B-c/2=-\lambda f(\ta),\qquad
\lambda\ \mathrm{constant},
\end{equation}
with $A:=(a^2+b^2)/(2a^2)$, $B:=(\al b-\bet a)/a^2$,
the metric
\[
\gK=-d(f(\ta)\kk^\flat)(J\cdot, \cdot)
\]
is K\"ahler-Einstein on the region $\{ -f\iota>0, f'>0 \}$.
\end{thm}
\begin{proof}
The complexified connection $1$-forms for the metric are given by
\Ref{gammas}, so we need only consider how
the computation of the Ricci form \Ref{rho} on our frame changes for
our current $f(\ta)$.

The Ricci form $\rho_{\scriptscriptstyle K}$ no longer annihilates
$\kk$, $\tT$ as the coefficients of $\hat\kk$, $\hat\tT$ in the connection
$1$-forms $\Gamma_1^1$, $\Gamma_2^2$ are no longer constant. We compute
the nonzero terms, using e.g. $(d\ta\we\hat \kk)(\kk,\tT)=-b$ etc.
\begin{align*}
\rho_{\scriptscriptstyle K}(\kk,\tT)&=i\left(\frac{f''}{2f'}\right)'((a-ib)(-b)+(b+ia)a)\\
&+i\left(\frac{f'}{2f}\right)'((a-ib)(-b)+(b+ia)a)\\
&=-(a^2+b^2)\left(\frac{f''}{2f'}+\frac{f'}{2f}\right)'\\
&=-\frac{a^2+b^2}2\left(\frac{(f^2)''}{(f^2)'}\right)'
\end{align*}
We now evaluate the Ricci form on $\xx$, $\yy$.
The exterior derivatives of terms involving $f(\ta)$ will only contribute terms of the form $d\ta\we\hat\kk$ or $d\ta\we\hat\tT$, which do not contribute when evaluated on $\xx$, $\yy$.
On these frame fields the only change in the calculation is to append the expressions involving
$f$ to the appropriate terms we have calculated in Theorem~\ref{central}. Thus
\begin{align*}
\rho_{\scriptscriptstyle K}(\xx,\yy)&=
\fr\iota{a^2}(a^2+b^2)\left(\frac{f''}{2f'}+\frac{f'}{2f}\right)+\iota\fr{-b\alpha+a\beta}{a^2}
-\fr12(d_\xx d_\xx+d_\yy d_\yy)(\log|\iota|)\\
&=\fr\iota{2a^2}(a^2+b^2)\frac{(f^2)''}{(f^2)'}+\iota\fr{-b\alpha+a\beta}{a^2}
-\fr12(d_\xx d_\xx+d_\yy d_\yy)(\log|\iota|)
\end{align*}
After verifying routinely that the Ricci form still vanishes for the current $f(\ta)$
on a pair of frame fields, one from $\VV$ and the other from $\HH$, we equate
the last displayed equation to $\lambda\gK(\xx,\xx)=-\lambda f(\ta)\iota$, and employ
\Ref{pde} to obtain \Ref{KEceq} after cancelling $\iota$ throughout the equation. On the other hand, equating our expression for $\rho_{\scriptscriptstyle K}(\kk,\tT)$ with $\lambda\gK(\kk,\kk)=\lambda a^2f'$ simply yields equality of the $\ta$-derivatives of both sides of \Ref{KEceq}. Our proof is complete, with the domain of the K\"ahler-Einstein metric obtained as usual by requiring positivity of $\gK(\kk,\kk)=a^2f'$ and $\gK(\xx,\xx)=-f\iota$.
\end{proof}

All the Lorentzian examples we had for central metrics also yield K\"ahler-Einstein metrics.
For instance, truncated pp-waves have $A=1$, $B=0$, so if $k(x,y)$, $h(x,y)$ are
chosen so that $\iota$ is a negative constant, or more generally so that $c=0$, then
the solution to the ODE \Ref{KEceq} is given implicitly as $F(f(\ta))=\ta+c_2$, where $c_2$
is constant and
\[ F(x)=\int\frac {3x}{-\lambda x^3+c_1}\,dx, \]
where $c_1$ is also a constant. Taking $c_1=c_2=0$ and $\lambda =-1$, this yields $f(\ta)=-3/\ta$,
and the K\"ahler metric is defined on $\{\ta<0\}$.
For the plane wave  $A=1/2$ and $B=0$, and recalling that $\ta=u$,
the same constants of integration and $\lambda$ choices
yield $f(u)=-3/(2u)$, with the domain $\{u<0\}$.
Returning to the truncated pp-wave case, if we choose $\iota$ to solve \Ref{pde}
for nonzero $c$, then the solution with vanishing constants of integration takes
the form \[ f(\ta)=\frac{3ce^{c\ta/4}}{4\lambda e^{c\ta/4}-1} \] in an appropriate
domain.

Turning to the Lie group examples of subsection~\ref{Lie-g}, those can
similarly be employed for K\"ahler-Einstein metrics, but the latter
will now be only {\em biconformal} to a left invariant metric when $\iota$ is constant:
\[ \gK=-f(\ta)\iota(\hat \xx^2+\hat \yy^2)+a^2f'(\ta)(\hat \kk^2+\hat \tT^2). \]

Computation of the sectional curvature of these K\"ahler metrics shows that
they are not all of constant holomorphic sectional curvature. More explicitly,
the solution of the ODE involves inverting an integral of a rational function.
The latter depends on integration constants. When these are set to
zero, the inversion is immediate, and the solution has constant holomorphic
sectional curvature, whereas otherwise it does not. More precisely,
one easily calculates from the covariant derivative expressions the following
sectional curvatures:
\begin{align*}
\text{$K_{\scriptscriptstyle K}(\xx,\yy)$}&\text{$=-2A\frac {f'(\ta)}{f^2(\ta)}+B\frac 1{f(\ta)}$\quad
for constant $\iota$,}\\
K_{\scriptscriptstyle K}(\kk,\tT)&=-A\left(\frac {f''(\ta)}{f'(\ta)}\right)'\fr 1{f'(\ta)}.\\
\end{align*}
But integrating equation~\Ref{KEceq} once and rearranging yields
\[
\text{$2A\frac {f'(\ta)}{f^2(\ta)}-B\frac 1{f(\ta)}=-\frac 23\lambda+c_1f^{-3}(\ta)$
for a constant $c_1$}.
\]
Comparing this equation and $K_{\scriptscriptstyle K}(\xx,\yy)$, it immediately
follows that a necessary condition for a K\"ahler-Einstein metric
obtained as in Theorem~\ref{KEc} with constant $\iota$
to have constant holomorphic sectional curvature is $c_1=0$. It is easy to see
this condition is also sufficient.

\subsection{Completeness and incompleteness of the K\"ahler-Einstein metrics}

We first determine in/completeness of the K\"ahler-Einstein metrics in the truncated pp-wave case,
with constant twist $\iota$. Our work is a minor elaboration on a method communicated to us
by Robert Ream.

Recall that for the truncated pp-waves, the frame is given by
\[
\kk = \partial_v + \partial_\ta,\quad
\tT = \partial_\ta,\quad
\xx = \partial_x + k\partial_v,\quad
\yy = \partial_y + h\partial_v
\]
for smooth functions $h = h(x, y)$, $k = k(x, y)$ satisfying $h_x - k_y = -1$.
The non-metric dual coframe is then given by
\[
\hat\kk = dv - kdx - hdy,\quad
\hat\tT = dv - kdx - hdy - d\ta,\quad
\hat\xx = dx,\quad
\hat\yy = dy.
\]

Let $\gamma(t) = (x(t), y(t), v(t), \ta(t))$ be a curve defined on an interval $I$ with left
endpoint $t_0$, and tangent given by
\begin{multline*} \gamma'= x'\partial_x + y'\partial_y + v'\partial_v + \ta'\partial_\ta
= \hat\xx(\gamma')\xx+ \hat\yy(\gamma')+\hat\kk(\gamma')\kk+ \hat\tT(\gamma')\tT\\
=x'\xx+y'\yy+(v'-kx'-hy')\kk+(v'-kx'-hy'-\ta')\tT.
\end{multline*}
As $a=1$ and we choose for simplicity $h$ and $k$ so that $\iota=-1$, the K\"ahler metrics
take the form
\be\lb{gK-KE} \gK=f(\ta)(\hat\xx^2+\hat\yy^2)+f'(\ta)(\hat\kk+\hat\tT)^2, \end{equation}
for a smooth function $f(\ta)$ considered in the domain $M=\{f>0, f'>0\}$.
By the Cauchy-Schwarz inequality, for any unit vector field $\mathbf u$,
\[ L(\gamma)=\int_I||\gamma'||\,dt\geq\int_I|\langle\gamma',\mathbf u\rangle|\,dt. \]
where $L(\gamma)$ is the length of $\gamma$. Applying this to unit vector fields in the direction of our
frame fields yields
\begin{align}
L(\gamma)&>\int_I\sqrt{f(\ta)}|x'|\,dt\geq \inf_{t\in I}\sqrt{f(\ta)}\left|\int_Ix'\,dt\right|,
\lb{11}\\
L(\gamma)&>\int_I\sqrt{f(\ta)}|y'|\,dt\geq \inf_{t\in I}\sqrt{f(\ta)}\left|\int_Iy'\,dt\right|,
\lb{22}\\
L(\gamma)&>\int_I\sqrt{f'(\ta)}|v'-kx'-hy'|\,dt\geq \left|\int_I\sqrt{f'(\ta)}(v'-kx'-hy')\,dt\right|=:L_1,
\lb{33}\\
L(\gamma)&>\int_I\sqrt{f'(\ta)}|v'-kx'-hy'-\ta'|\,dt\geq \left|\int_I\sqrt{f'(\ta)}(v'-kx'-hy'-\ta')\,dt\right|\nonumber\\
&=\left|\pm L_1+\int_I\sqrt{f'(\ta)}\ta'\,dt\right|,\lb{44}
\end{align}
Suppose for a given initial condition, a solution $f(\ta)$ to \Ref{KEceq} with $A=1$, $B=0$, $f(\ta)>0$
and $f'(\ta)>0$ is defined on a (possibly unbounded) maximal interval $(\ta_a,\ta_b)$,
and for any $\ta_0$ in this interval,
\be\lb{length} \text{ $\left|\int_{\ta_0}^{\ta_f}\sqrt{f'(\ta)}\,d\ta\right|=\infty$, }\qquad
\text{ for both $\ta_f=\ta_a$, $\ta_f=\ta_b$.} \end{equation}
Now let $\ta_0=\ta(t_0)$, $\ta_p=\lim_{t\to\sup I}\ta(t)$.
If the curve length $L(\gamma)$ is finite, then \Ref{44} and \Ref{length} show that $\ta_p$ is different
from the boundary values $\ta_a$, $\ta_b$, so that $\ta(t)$, $t\in I$ is bounded away
from $\ta_a$, $\ta_b$. Then inequalities \Ref{11}, \Ref{22} show that $x(t)$, $y(t)$ are bounded,
hence so are $k(x(t),y(t))$, $h(x(t),y(t))$, and by \Ref{33}, also $v(t)$.

Now assume that the curve $\gamma$ is escaping, i.e. it leaves every compact set of $M$.
Then one of the following occur: $\ta(t)$ approaches either $\ta_a$ or $\ta_b$, $x(t)\to\pm\infty$,
$y(t)\to\pm\infty$, or $v(t)\to\pm\infty$. By the preceding paragraph, this can only happen
when $L(\gamma)=\infty$, proving that $\gK$ is complete.

On the other hand, if one of the integrals in \Ref{length} are finite, then
the curve $\gamma(\ta)=(0,0,0,\ta)$, parametrized by $\ta$  between some $\ta_0$
and the appropriate $\ta_f$  has length given by
$\int_{\ta_0}^{\ta_f}\sqrt{f'(\ta)}\,d\ta<\infty$, hence $\gK$ is incomplete.

We examine the state of affairs for some of our K\"ahler-Einstein metrics associated
to a truncated pp-wave. The explicit solution $f(\ta)=-3/\ta$ to \Ref{KEceq}, with $\lambda=-1$, is
positive with a positive derivative on $(\ta_a,\ta_b)=(-\infty,0)$, and the integrals in \Ref{length}
are both infinite in this case. Hence $\gK$ is complete. Together with the
fact that we have seen it is of constant holomorphic sectional curvature, this shows
$\gK$ is isometric to the standard metric on complex hyperbolic space with holomorphic
sectional curvature $-2/3$.

Now suppose $f(\ta)$ solves \Ref{KEceq} and satisfies
\be\lb{complt} f'=-\frac{(2\lambda/3) f^3+c_1}{2f}, \end{equation}
for constants $\lambda<0$, $c_1>0$. Then whenever $f>0$, we also have
$f'>0$ so long as
\be\lb{bound}
f>[-3c_1/2\lam]^{1/3}>0.
\end{equation}
Moreover, $f''=(f'/4f^2)(-8\lam f^3/3+2c_1)>0$ as well so that $f'$ is also increasing.
Suppose such a solution, satisfying \Ref{bound}, with an initial condition consistent with
this inequality,  is maximally defined on $(\ta_a,\ta_b)$.
Then for $\ta_0\in (\ta_a,\ta_b)$, with $f_0=f(\ta_0)$ and $f_a=\lim_{\ta\searrow\ta_a}f(\ta)$, we have
\begin{align*}
\left|\int_{\ta_0}^{\ta_a}\right.&\left.\sqrt{f'(\ta)}\,d\ta\right|=\\
&\left|\int_{f_0}^{f_a}\sqrt{\frac 1{\ta'(f)}}\ta'(f)\,df\right|=\left|\int_{f_0}^{f_a}\sqrt{\ta'(f)}\,df\right|=\\
&\left|\int_{f_0}^{f_a}\sqrt{-\frac{2f}{(2\lambda/3) f^3+c_1}}\, df\right|<\infty
\end{align*}
since inside the square-root in the integrand, the rational function has a simple pole
at $f_q:=[-3c_1/2\lam]^{1/3}$, and $f_q\le f_a<\infty$. Thus the corresponding metric
$\gK$ is incomplete.

For the purpose of the material in the rest of this subsection
we mention that in this setting $\ta_b$ must be finite and
$f_b:=\lim_{\ta\nearrow\ta_b}f(\ta)=\infty$. Thus $\sqrt{\ta'(f)}$
vanishes asymptotically as $1/f$ as $f$ approaches $f_b$ and
its integral from $f_0$ to $f_b$ diverges. Hence completeness
is unhindered near $\ta_b$. Addtionally, for an autonomous ODE
of the form \Ref{complt} a maximal solution with an initial condition
$f(\tau_0)>0$ is known to extend to $\ta=-\infty$, with $f$ approaching $f_q$,
so that we can assume $f_a=f_q$.

However, we now demonstrate that it is possible to obtain complete K\"ahler-Einstein
metrics for $f$ satisfying \Ref{complt} for some values of $c_1\ne 0$, via a method due
again to Robert Ream, relying on recent results in~\cite{vz}. The metrics will be of
cohomogeneity one, for the action of a group $\mathcal{G}=\mathrm{nil_3}/\tilde{\mathbb{Z}}$, which
is the quotient of the Heisenberg group  by an infinite cyclic group lying in its center,
given by
\[ \tilde{\mathbb{Z}}:= \left\{\begin{bmatrix} 1 & 0 & \ell n\\ 0 & 1 & 0\\ 0 & 0 & 1  \end{bmatrix} |\ n\in\mathbb{Z} \right\},\]
where the choice of a real number $\ell$ is determined by later considerations.
$\mathcal{G}$ has center $K$ isomorphic to
$SO(2)$, whose transitive action on the circle $S^1$ extends to a linear action on $V:=\mathbb{R}^2$. We consider the homogeneous vector bundle $M=\mathcal{G}\times_K V$ (in which points of the product
are identified according to $(g,v)\sim(gk^{-1},kv)$ for $k\in K$).
$\mathcal{G}$ acts on $M$ by left multiplication on the first factor.
The action of $\mathcal{G}$ has trivial isotropy at points of a regular orbit, but isotropy $K$ at a
point of the singular orbit $\mathcal{G}/K\approx \mathbb{R}^2$.

A left-invariant frame for $\mathcal{G}$ has the form $\xx=\partial_x$, $\yy=\partial_y+x\partial_z$,
$\tT=\partial_z$, to which we will add $\kk=\partial_\ta$. The corresponding left invariant coframe
is $\hat\kk=d\ta$, $\hat\tT=dz-xdy$, $\hat\xx=dx$, $\hat\yy=dy$.  We define the metric
$\gK$ by the formula \Ref{gK-KE} for $f$ satisfying \Ref{complt}. We also require as before
$\lambda<0$, $c_1>0$ and $f>(-3c_1/2\lambda)^{1/3}$.

Changing $\ta$ to the variable
$f$ and noting that $f'(\tau)(d\tau^2+\hat{t}^2)=\frac{df^2}{f'(\tau)}+f'(\tau)\hat{t}^2$,
this metric can be written in the form
\begin{multline*}
\gK=\frac{2f}{-\frac{2}{3\lambda}f^3-c_1}df^2+\frac{-\frac{2}{3\lambda}f^3-c_1}{2f}
\hat\tT^2+f(\hat{\xx}^2+\hat{\yy}^2)\\
=\frac{2f}{p(f^3-\alpha^3)}df^2+\frac{p(f^3-
\alpha^3)}{2f}\hat{\tT}^2+ f(\hat{\xx}^2+\hat{\yy}^2)
\end{multline*}
defined on the domain $f\in (\al,\infty)$, where $\al=(-3c_1/2\lambda)^{1/3}$ and $p=-2/3\lambda$.

For a curve $\gamma$ as before, we have $\gK(\gamma',\xx/||\xx||)=x'\sqrt{f}$, $\gK(\gamma',\yy/||\yy||)=y'\sqrt{f}$, $\gK(\gamma',\tT/||\tT||)=(z'-xy')\sqrt{f'}$,
$\gK(\gamma',\kk/||\kk||)=\ta'\sqrt{f'}$, so that the argument for reducing completeness
to the validity of \Ref{length} would still hold. Note that taking the quotient group
makes one of the coordinates automatically bounded, so that argument works even more
easily.

Of course we have seen that \Ref{length} does not hold at $\ta_a$. This difficulty
could be overcome, provided that we can attach smoothly a singular orbit at $f=\alpha$,
so that a curve for which $f$ approaches $\ta_a$ with the other coordinates bounded
will no longer be an escaping curve.

Smoothness conditions for a metric at a singular orbit have been given recently
in \cite{vz}. To apply these we need to write the metric near $\alpha$ in the form
$dt^2+h_t$ where $t=0$ corresponds to $f=\al$. Asymptotically near $\al$ we have
$dt=\sqrt{\frac{2f}{p(f^3-\al^3)}}df\approx \sqrt{\frac{2}{3p\al(f-\al)}}df$ so
that $t\approx\sqrt{q(f-\al)}$ for $q=8/(3p\al)$.

Thus near $\al$
\[
\gK\approx dt^2+\frac{3p\al^2t^2/q}{2(t^2/q+\al)}\hat\tT^2+\frac 1q(t^2+q\alpha)(\hat\xx^2+\hat\yy^2).
\]

We compare this with the smoothness conditions in \cite{vz}, which in our case, for $\mathfrak{m}=\mathrm{span}(\xx,\yy)$, $\mathfrak{p}=\mathrm{span}(\tT)$,
are, near $t=0$,
\begin{align*}
&\text{$\gK(\mathfrak{m},\mathfrak{m})$ is even in $t$,}\\
&\text{$\gK(\mathfrak{p},\mathfrak{m})=t^2\phi(t^2)$,}\\
&\text{$\gK(X,X)=t^2+t^4\psi(t^2)$ for $X\in\mathfrak{p}$ such that $\gK(X,X)|_{t=0}=1$.}
\end{align*}
Only the last condition is not automatic in our case. The coefficient of $\hat\tT^2$ is
\[ \frac{3p\al^2t^2/q}{2(t^2/q+\al)}=\frac{3p\al t^2}{2q}(1-t^2/\al q+O(t^4)). \]
Since the part in brackets is also even, this coefficient will have the right form if $\frac{3p\al}{2q}=1$, i.e. $c_1=(2\lambda)^4/3$.

The K\"ahler form similarly extends smoothly to the singular fiber. In fact, it is
\begin{multline*} df\we\hat t+ f(\xx\we\yy)=
d\left[\left(\sqrt{f-\al}\right)^2\right]\we\hat \tT+ f(\hat\xx\we\hat\yy)\\
\approx q^{-1}d(t^2)\we\hat \tT+ q^{-1}(t^2+\al)(\hat\xx\we\hat\yy),
\end{multline*}
whereas modifying the conditions in \cite{vz} so that they apply to a $2$-form, shows that
in our case smoothness requires that near $t=0$ the coefficient of $dt\we\hat \tT$ has
the form $t\phi(t^2)$ and the coefficient of $\hat\xx\we\hat\yy$ is even. Thus $\gK$
is a complete K\"ahler-Einstein metric.

Some of these methods can be applied to our other Lie groups. They yield, for example, the
complete biaxial K\"ahler-Einstein metrics of \cite{g-p}, \cite{pe} and \cite{d-s1}, which are
cohomogeneity one for $SU(2)$. For the remaining Lie groups the methods based on \cite{vz}
can still be applied, but estimates along the lines of \Ref{11}-\Ref{44} do not lead to the
same level of control of the coordinates of a finite length curve. We hope to address
this in future work.

\section{K\"ahler-Einstein metrics associated to Lorentzian warped products}\lb{sec:3}

In this section we assume $M=N\times\mathbb{R}$, where the real line
is equipped with a coordinate function $\ta$ and $N$ is a $3$-manifold.
Furthermore, $N$ admits a Riemannian metric $\bar{g}$ with
\be\lb{vfield}
\begin{aligned}
&\text{a unit length vector field $\bar{\kk}$, which is geodesic,}\\
&\text{shear-free, and has a nowhere vanishing twist function $\bar{\iota}$.}
\end{aligned}
\end{equation}
Here the shear operator and twist
function are defined for a vector field on a $3$-manifold in complete analogy
with the four dimensional case (see \cite[subsection 2.4]{a-m}). For example,
the twist function of $\bar\kk$ is given by $\bar\iota=\bar{g}(\bar\kk,[\bar \xx,\bar\yy])$,
for an oriented $\bar g$-orthonormal frame $\bar\xx$, $\bar\yy$ of the
$\bar g$-orthogonal complement of $\bar\kk$.

We now list some of the consequences of \cite[Sec.~8]{a-m} that pertain to
admissible Lorentzian warped product metrics on $M$.
By 
Theorem~$5$ in \cite{a-m} and its proof, after choosing an appropriate orientation on $M$
and a positive $C^\infty$ function $w(\ta)$ on $\mathbb{R}$,
the metric $g=-d\ta^2+w(\ta)^2\bar{g}$ is admissible, with $J=J_{g,\kk,\tT}$ sending  the null vector
field $\kk=\bar{\kk}/w+\partial_\ta$, which is geodesic or strictly pre-geodesic, to the gradient field
$\tT=\nabla\ta=-\partial_\ta$. Our interest will be in the case where $w$ is nonconstant, in which $\kk$ is actually strictly pre-geodesic: $\n_\kk\kk=(w'/w)\kk$ \cite[formula~(32)]{a-m}.

As $g$ is admissible, it gives rise to K\"ahler metrics of the form
$\gK=-d(f(\ta)\kk^\flat)(J\cdot,\cdot)$, for smooth $f$ defined on $\mathbb{R}$,
in any region where the inequalities
\be\lb{region}
f\iota<0,\qquad f'+fw'/w>0
\end{equation}
hold, where $\iota$ is, as usual, the twist function of $\kk$, and the prime denotes $d/d\ta$.
Note that \cite[Theorem 5]{a-m} only pertains to the special case $f(\ta)=e^\ta$,
but \Ref{region} appears in essence in \cite[Sec.~8.5, formula~(41)]{a-m}.

The next three facts appear in \cite[subsection~8.1]{a-m}. The relation between the two twists is
\be\lb{twists}
\iota=w^{-1}\bar\iota,
\end{equation}
Similarly, the vanishing of the $\bar g$-shear of $\bar\kk$ implies
\be\lb{shearAM}
\text{the $g$-shear of $\kk$ vanishes.}
\end{equation}
Finally, the Lorentzian metric values yield
\be\lb{Lval}
\text{$d_\kk\ta=g(\kk,\tT)=1$, $d_\tT\ta=g(\tT,\tT)=-1$ while clearly $d_\xx\ta=d_\yy\ta=0$.}
\end{equation}

We assume $N$ admits a system of local $\bar{g}$-orthonormal frames of the form $\bar{\xx}$, $\bar{\yy}$  for the $\bar{g}$-orthogonal complement of $\bar{\kk}$, such that for each one of them
\be\lb{N-cond}
\begin{aligned}
&\text{$\bar{\iota}$ is negative,}\\
&\text{$[\bar{\kk},\bar{\xx}]=\al\bar{\yy},\quad [\bar{\kk},\bar{\yy}]=-\al\bar{\xx},\quad $ for a (frame independent) constant $\al$,}\\
&\text{$[\bar{\xx},\bar{\yy}]$ is a multiple of $\bar{\kk}$, hence necessarily $\bar{\iota}\bar{\kk}$.}
\end{aligned}
\end{equation}
If $\bar\iota$ is constant, $\bar{\kk}$, $\bar{\xx}$, $\bar{\yy}$ form a $3$-dimensional Lie algebra
with respect to the Lie bracket, and in principle one can consider others (cf. \cite{bow} and the
appendix of \cite{derdi}). The corresponding connected simply connected Lie group in that case
is $S^3$ (if $\al\ne 0$) or $\mathbb{R}^3$ (if $\al=0$). These $3$-manifolds will be
referred to in our examples.

Even if $\bar\iota$ is not constant, we always have
\[
d_{\bar \kk}\bar\iota=0.
\]
This follows from the Newman-Penrose related equation $d_{\bar \kk}\bar\iota=
-(\bar\delta\bar \kk)\bar\iota$, valid for the unit length geodesic field $\bar{\kk}$,
since the divergence is just $\bar\delta\bar \kk=\bar g(\n_{\bar \xx}\bar \kk,\bar \xx)+
\bar g(\n_{\bar \yy}\bar \kk,\bar \yy)+\bar g(\n_{\bar \kk}\bar \kk,\bar \kk)=
\bar g([\bar \xx,\bar \kk],\bar \xx)+\bar g([\bar \yy,\bar \kk],\bar \yy)=0$ by \Ref{N-cond}.
This calculation also shows that under assumptions \Ref{vfield} and \Ref{N-cond},
$\bar{\kk}$ is in fact a Killing field (see \cite[Lemma 2.6]{a-m}).

We can now state our theorem.
\begin{thm}\lb{ke}
Let $g$ be an admissible warped product metric, constructed as above from a
warping function $w(\ta)$ and a Riemannian $3$-manifold $(N,\bar{g})$ satisfying
\Ref{vfield} and \Ref{N-cond}. Also let $f(\tau)$ be a real valued function on
$\mathbb{R}$ and $I\subset\mathbb{R}$ be given by the inequalities
\be\lb{ineq}
 f>0,\qquad (fw)' > 0.
\end{equation}
Then $g$ and $f(\tau)$ induce an associated K\"ahler-Einstein
metric\linebreak $\gK=-d(f(\ta)\kk^\flat)(J\cdot,\cdot)$ with Einstein constant $\lambda$ on
$N\times I$ if and only if there exists a constant $C$ such that
\be\lb{KE-N}
\text{$\left(d_{\bar\xx}^2+d_{\bar\yy}^2\right)\log|\bar\iota|=-2\lambda C\bar\iota$\quad
holds on $N$}
\end{equation}
and
\be\lb{KE}
\text{$\fr{(fw)''}{(fw)'}+2\fr{w'}w+\fr{f'}{f}+\fr{\al}w=-\lambda \left(\frac Cw+f\right)$\quad
holds on $I$.}
\end{equation}
\end{thm}
Inequalities \Ref{ineq} are just \Ref{region}, as is clear from \Ref{twists},
$w>0$ and $\bar\iota<0$.

The proof will run similarly to that of Theorem~\ref{central}, and
we break it again into subsections.

\subsection{Basic relations and the connection}\lb{basic}
We first lift $\bar{\xx}$, $\bar{\yy}$ to $M$, giving a $g$-orthonormal basis
$\xx=\bar{\xx}/w$, $\yy=\bar{\yy}/w$ of $\HH=\mathrm{span}(\kk,\tT)^\perp$. Then
we verify most of the $2$-dependency described in the introduction, with respect to the functions
$\ta$ and $\bar\iota$, except that $d_\xx\bar\iota$ and $d_\yy\bar\iota$ remain unspecified.
On the other hand clearly $d_\tT\bar\iota=\partial_\ta\bar\iota=0$ and $d_\kk\bar{\iota}=w^{-1}d_{\bar \kk}\bar\iota-d_\tT\bar\iota=0$. The directional derivatives of $\ta$ are computed in \Ref{Lval}.
Turning to frame Lie brackets,
\[
[\kk,\xx]=[\bar{\kk}/w+\partial_\ta,\bar{\xx}/w]=w^{-2}[\bar{\kk},\bar{\xx}]-w^{-2}w'\bar{\xx}
=\al w^{-1}\yy-w^{-1}w'\xx,
\]
as $d_{\bar{\xx}}w^{-1}=d_{\bar{\kk}}w^{-1}=0$ and $[\partial_\ta,\bar{\xx}]=0$.
Similarly
\[
[\kk,\yy]=-\al w^{-1}\xx-w^{-1}w'\yy,\qquad [\tT,\xx]=w^{-1}w'\xx,\qquad [\tT,\yy]=w^{-1}w'\yy.
\]
Finally
\[
\begin{aligned}
&[\xx,\yy]=[\bar{\xx}/w,\bar{\yy}/w]=
[\bar{\xx},\bar{\yy}]/w^2=\bar{\iota}\bar{\kk}/w^2=\bar{\iota}w^{-1}(\kk+\tT),\\
&[\kk,\tT]=[\bar{\kk}/w+\partial_\ta,-\partial_\ta]=[\partial_\ta,\bar{\kk}/w]
=-w^{-2}w'\bar{\kk}=-w^{-1}w'(\kk+\tT).
\end{aligned}
\]

Next, the K\"ahler metric values on our $\gK$-orthogonal basis are
\[
\begin{aligned}
&\gK(\xx,\xx)=\gK(\yy,\yy)=-f\iota=-fw^{-1}\bar{\iota},\\
&\gK(\kk,\kk)=\gK(\tT,\tT):=c=-f'G/\ell+fw^{-1}w' g(\kk,\tT)=f'+fw^{-1}w'=(fw)'/w.
\end{aligned}
\]
The first line follows from the last part of \Ref{gK} and \Ref{twists}.
The second line results from the second line of \Ref{gkk}, with $\gamma=w'/w$,
using \Ref{Lval} which yields $G=-1$, while $\ell=1$ as $\tT=\n\ta$.

We wish to employ \Ref{Kosz} and the Koszul formula. We have, for example
$2\gK(\nK_\kk\kk,\kk)=d_\kk\gK(\kk,\kk)=c'$, $2\gK(\nK_\kk\kk,\tT)=-d_\tT\gK(\kk,\kk)-2\gK(\kk,[\kk,\tT])
=c'+2w^{-1}w'c$, $2\gK(\nK_\kk\kk,\xx)=2\gK(\nK_\kk\kk,\yy)=0$. This,
along with the connection commuting with the complex structure, and
the value of $[\kk,\tT]$ yields the four relations
\[
\begin{aligned}
\nK_\kk\kk&=\fr{c'}{2c}(\kk+\tT)+w^{-1}w'\tT,\qquad
\nK_\tT\tT=-\fr{c'}{2c}(\kk+\tT)-w^{-1}w'\kk,\\
\nK_\kk\tT&=\fr{c'}{2c}(\tT-\kk)-w^{-1}w'\kk,\qquad
\nK_\tT\kk=\fr{c'}{2c}(\tT-\kk)+w^{-1}w'\tT.
\end{aligned}
\]
We note here that $c'/c=(fw)''/(fw)'-w'/w$.

Given that $d_\tT\bar\iota=d_\kk\bar\iota=0$, we have
\[
\begin{aligned}
d_\kk\gK(\xx,\xx)&=d_\kk\gK(\yy,\yy)=-d_\tT\gK(\xx,\xx)=-d_\tT\gK(\yy,\yy)\\[4pt]
&=-(f\iota)'=-f'\iota-f(w^{-1}\bar{\iota})'
=-f'\iota+fw^{-2}w'\bar{\iota}.
\end{aligned}
\]
Using notations of subsection~\ref{conn-ct}, we note that
\[
\iota_{\scriptscriptstyle K}=\iota_{\scriptscriptstyle K}^\tT=
\bar{\iota}w^{-1}c=\bar{\iota}w^{-2}(fw)'=\bar{\iota}w^{-1}f'+\bar{\iota}w^{-2}w'f
=f'\iota+\bar{\iota}w^{-2}w'f
\]

We have
$2\gK(\nK_\xx\kk,\xx)=d_\kk\gK(\xx,\xx)-2\gK(\xx,[\kk,\xx])$,
$2\gK(\nK_\xx\kk,\yy)=-\iota_{\scriptscriptstyle K}$,
$2\gK(\nK_\xx\kk,\kk)=0$, $2\gK(\nK_\xx\kk,\tT)=-\gK(\xx,[\kk,\tT])=0$,
where for the second equality we used \Ref{shearAM} and Remark~\ref{shearK}
as in the previous section.
Deducing from this also the expressions for $\nK_\yy\kk$, $\nK_\xx\tT$, $\nK_\yy\tT$, we have
\[
\begin{aligned}
\nK_\xx\kk&=\left(\fr{f'}{2f}+\fr{w'}{2w}\right)(\xx+\yy),
\qquad
\nK_\yy\kk=\left(\fr{f'}{2f}+\fr{w'}{2w}\right)(\yy-\xx),\\
\nK_\xx\tT&=\left(\fr{f'}{2f}+\fr{w'}{2w}\right)(\yy-\xx),
\qquad
\nK_\yy\tT=\left(\fr {f'}{2f}+\fr{w'}{2w}\right)(-\xx-\yy).
\end{aligned}
\]

In the same mode we have, for example
$2\gK(\nK_\xx\xx,\kk)=-d_\kk\gK(\xx,\xx)-2\gK(\xx,[\xx,\kk])$,
$2\gK(\nK_\xx\xx,\tT)=-d_\tT\gK(\xx,\xx)-2\gK(\xx,[\xx,\tT])$,
$2\gK(\nK_\xx\xx,\xx)=d_\xx\gK(\xx,\xx)$, $2\gK(\nK_\xx\xx,\yy)=-d_\yy\gK(\xx,\xx)-2\gK(\xx,[\xx,\yy])=0$.
Hence we compute
\[
\begin{aligned}
&\nK_\xx\xx=\frac{f'\iota+fw^{-2}w'\bar{\iota}}{2c}(\kk-\tT)
+\frac{d_\xx\bar\iota}{2\bar\iota}\xx-\frac{d_\yy\bar\iota}{2\bar\iota}\yy,\\
&\nK_\yy\yy=\frac{f'\iota+fw^{-2}w'\bar{\iota}}{2c}(\kk-\tT)
-\frac{d_\xx\bar\iota}{2\bar\iota}\xx+\frac{d_\yy\bar\iota}{2\bar\iota}\yy.\\
&\nK_\xx\yy=\frac{f'\iota+\bar{\iota}w^{-2}w'f}{2c}(\kk+\tT)
+\frac{d_\yy\bar\iota}{2\bar\iota}\xx+\frac{d_\xx\bar\iota}{2\bar\iota}\yy,\\
&\nK_\yy\xx=-\frac{f'\iota+\bar{\iota}w^{-2}w'f}{2c}(\kk+\tT)
+\frac{d_\yy\bar\iota}{2\bar\iota}\xx+\frac{d_\xx\bar\iota}{2\bar\iota}\yy,
\end{aligned}
\]
Finally, using formulas like $\nK_\kk\xx=\nK_\xx\kk+[\kk,\xx]$ we have
\[
\begin{aligned}
\nK_\kk\xx&=\nK_\xx\kk+\al w^{-1}\yy-w^{-1}w'\xx,
\qquad
&&\nK_\kk\yy=\nK_\yy\kk-\al w^{-1}\xx-w^{-1}w'\yy,\\
\nK_\tT\xx&=\nK_\xx\tT+w^{-1}w'\xx,
\qquad
&&\nK_\tT\yy=\nK_\yy\tT+w^{-1}w'\yy.
\end{aligned}
\]
From these calculations we can write the complex-valued connection $1$-forms
as in the previous section:
\[
\begin{aligned}
\Gamma_1^1&=\left(\fr{c'}{2c}(1+i)+i\fr{w'}w\right)\hat{\kk}+
\left(\fr{c'}{2c}(i-1)+i\fr{w'}w\right)\hat{\tT},\\
\Gamma_1^2&=\left(\fr{f'}{2f}+\fr{w'}{2w}\right)((1+i)\hat{\xx}+(i-1)\hat{\yy}),\\
\Gamma_2^1&=\fr{f'\iota+fw^{-2}w'\bar{\iota}}{2c}((1-i)\hat{\xx}-(1+i)\hat{\yy}),\\
\Gamma_2^2&=\left[\left(\fr{f'}{2f}+\fr{w'}{2w}\right)(1+i)+i\fr{\al}w-\fr{w'}w\right]\hat{\kk}
+\left[\left(\fr{f'}{2f}+\fr{w'}{2w}\right)(i-1)+\fr{w'}w\right]\hat{\tT}\\
&+\frac 12\left(d_\xx\log|\bar\iota|-id_\yy\log|\bar\iota|\right)\hat\xx
+\frac 12\left(d_\yy\log|\bar\iota|+id_\xx\log|\bar\iota|\right)\hat\yy.
\end{aligned}
\]

\subsection{The Ricci form}
Since the coefficients of $\Gamma_i^j$ are functions of $\ta$, the Ricci form is
\be\lb{Ricci-w}
\begin{aligned}
\rho_{\scriptscriptstyle K}&=i(d\Gamma_1^1+d\Gamma_2^2)\\
&=i\left(\fr{c'}{2c}(1+i)+i\fr{w'}w\right)d\hat{\kk}+
i\left(\fr{c'}{2c}(i-1)+i\fr{w'}w\right)d\hat{\tT}\\
&+i\left[\left(\fr{f'}{2f}+\fr{w'}{2w}\right)(1+i)+i\fr{\al}w-\fr{w'}w\right]d\hat{\kk}\\
&+i\left[\left(\fr{f'}{2f}+\fr{w'}{2w}\right)(i-1)+\fr{w'}w\right]d\hat{\tT}\\
&+i\left(\fr{c'}{2c}(1+i)+i\fr{w'}w\right)'d\ta\we\hat{\kk}+
i\left(\fr{c'}{2c}(i-1)+i\fr{w'}w\right)'d\ta\we\hat{\tT}\\
&+i\left[\left(\fr{f'}{2f}+\fr{w'}{2w}\right)(1+i)+i\fr{\al}w-\fr{w'}w\right]'d\ta\we\hat{\kk}\\
&+i\left[\left(\fr{f'}{2f}+\fr{w'}{2w}\right)(i-1)+\fr{w'}w\right]'d\ta\we\hat{\tT}\\
&+\frac i2\left(d_\xx\log|\bar\iota|-id_\yy\log|\bar\iota|\right)d\hat\xx
+\frac i2\left(d_\yy\log|\bar\iota|+id_\xx\log|\bar\iota|\right)d\hat\yy\\
&+\frac i2d\left(d_\xx\log|\bar\iota|-id_\yy\log|\bar\iota|\right)\we\hat\xx
+\frac i2d\left(d_\yy\log|\bar\iota|+id_\xx\log|\bar\iota|\right)\we\hat\yy.
\end{aligned}
\end{equation}
To proceed further, we compute
\[
\begin{aligned}
&d\hat{\kk}(\kk,\tT)=-\hat{\kk}([\kk,\tT])=-\hat{\kk}(-w^{-1}w'(\kk+\tT))=w^{-1}w',\qquad
d\hat{\tT}(\kk,\tT)=w^{-1}w',\\
&d\hat{\kk}(\xx,\yy)=-\hat{\kk}([\xx,\yy])=-\hat{\kk}(\bar{\iota}w^{-1}(\kk+\tT))=
-\bar{\iota}w^{-1},\qquad
d\hat{\tT}(\xx,\yy)=-\bar{\iota}w^{-1},\\
&(d\ta\we\hat{\kk})(\kk,\tT)=1,\  (d\ta\we\hat{\tT})(\kk,\tT)=1,\
(d\ta\we\hat{\kk})(\xx,\yy)=0,\
(d\ta\we\hat{\tT})(\xx,\yy)=0,\\
&d\hat{\xx}(\kk,\tT)=-\hat{\xx}([\kk,\tT])=0,\
d\hat{\yy}(\kk,\tT)=0,\
d\hat{\xx}(\xx,\yy)=0,\
d\hat{\yy}(\xx,\yy)=0.
\end{aligned}
\]
Thus
\[
\begin{aligned}
\rho_{\scriptscriptstyle K}(\kk,\tT)&=i\left(\fr{c'}{2c}(1+i)+i\fr{w'}w\right)\fr{w'}w
+i\left(\fr{c'}{2c}(i-1)+i\fr{w'}w\right)\fr{w'}w\\
&+i\left[\left(\fr{f'}{2f}+\fr{w'}{2w}\right)(1+i)+i\fr{\al}w-\fr{w'}w\right]\fr{w'}w\\
&+i\left[\left(\fr{f'}{2f}+\fr{w'}{2w}\right)(i-1)+\fr{w'}w\right]\fr{w'}w\\
&+i\left(\fr{c'}{2c}(1+i)+i\fr{w'}w\right)'
+i\left(\fr{c'}{2c}(i-1)+i\fr{w'}w\right)'\\
&+i\left[\left(\fr{f'}{2f}+\fr{w'}{2w}\right)(1+i)+i\fr{\al}w-\fr{w'}w\right]'\\
&+i\left[\left(\fr{f'}{2f}+\fr{w'}{2w}\right)(i-1)+\fr{w'}w\right]'\\
\end{aligned}
\]
so that
\[
\begin{aligned}
\rho_{\scriptscriptstyle K}(\kk,\tT)
&=\left[\fr{c'}{2c}i2i+2i^2\fr{w'}w+\left(\fr{f'}{2f}+\fr{w'}{2w}\right)i2i+i^2\fr{\al}w\right]\fr{w'}w\\
&+\left(\fr{c'}{2c}\right)'i2i+2i^2\left(\fr{w'}w\right)'
+\left(\fr{f'}{2f}+\fr{w'}{2w}\right)'i2i+i^2\left(\fr{\al}w\right)'\\
&=\fr1w\left[\left[-\fr{c'}{c}-2\fr{w'}w-
\left(\fr{f'}{f}+\fr{w'}{w}\right)-\fr{\al}w\right]w\right]'\\
&=\fr1w\left[\left[-\fr{(fw)''}{(fw)'}-\fr{w'}w-
\left(\fr{f'}{f}+\fr{w'}{w}\right)-\fr{\al}w\right]w\right]'\\
&=-\fr1w\left[\left[\fr{(fw)''}{(fw)'}+2\fr{w'}w+
\fr{f'}{f}+\fr{\al}w\right]w\right]',
\end{aligned}
\]
where in the penultimate step we used $c'/c=(fw)''/(fw)'-w'/w$.

Similarly,
\[
\begin{aligned}
\rho_{\scriptscriptstyle K}(\xx,\yy)
&=i\left(\fr{c'}{2c}(1+i)+i\fr{w'}w\right)\left(-\fr{\bar{\iota}}w\right)
+i\left(\fr{c'}{2c}(i-1)+i\fr{w'}w\right)\left(-\fr{\bar{\iota}}w\right)\\
&+i\left[\left(\fr{f'}{2f}+\fr{w'}{2w}\right)(1+i)+i\fr{\al}w-\fr{w'}w\right]\left(-\fr{\bar{\iota}}w\right)\\
&+i\left[\left(\fr{f'}{2f}+\fr{w'}{2w}\right)(i-1)+\fr{w'}w\right]\left(-\fr{\bar{\iota}}w\right)\\
&+\frac i2\left(-d_\yy d_\xx\log|\bar\iota|+id_\yy d_\yy\log|\bar\iota|\right)
+\frac i2\left(d_\xx d_\yy\log|\bar\iota|+id_\xx d_\xx\log|\bar\iota|\right)\\
&=\left[-\fr{c'}{c}-2\fr{w'}w-
\left(\fr{f'}{f}+\fr{w'}{w}\right)-\fr{\al}w\right]\left(-\fr{\bar{\iota}}w\right)
-\frac12\left(d_\xx d_\xx+d_\yy d_\yy\right)\log|\bar\iota|\\
&=\left[\fr{(fw)''}{(fw)'}+2\fr{w'}w+
\fr{f'}{f}+\fr{\al}w\right]\fr{\bar{\iota}}w-\frac12\left(d_\xx d_\xx+d_\yy d_\yy\right)\log|\bar\iota|,
\end{aligned}
\]
where we have used the fact that $d_{[\xx,\yy]}\log|\bar\iota|=\bar\iota w^{-1}d_{\kk+\tT}\bar\iota=0$.

Note also that $\rho_{\scriptscriptstyle K}(\HH,\VV)=0$, since $d\hat{\kk}$, $d\hat{\tT}$,
$d\tau\we\hat{\kk}$, $d\ta\we\hat{\tT}$ all vanish on a pair of fields one from each of these distributions, whereas cancellations occur for the terms involving $\hat\xx$, $\hat\yy$ as in
the proof of Theorem~\ref{central}.

The last two computations for $\rho_{\scriptscriptstyle K}$ are of course just the Ricci
curvature values $\mathrm{Ric}_{\scriptscriptstyle K}(\xx,\xx)=
\mathrm{Ric}_{\scriptscriptstyle K}(\yy,\yy)$ and
$\mathrm{Ric}_{\scriptscriptstyle K}(\kk,\kk)=\mathrm{Ric}_{\scriptscriptstyle K}(\tT,\tT)$,
which we now compare with $\gK(\xx,\xx)=\gK(\yy,\yy)$ and $\gK(\kk,\kk)=\gK(\tT,\tT)$. Letting
$L=\fr{(fw)''}{(fw)'}+2\fr{w'}w+\fr{f'}{f}+\fr{\al}w$, the $\xx, \xx$ equation easily yields
$L=\frac w{2\bar\iota}\left(d_\xx^2+d_\yy^2\right)\log|\bar\iota|-\lambda f$, which,
substituted into the $\kk, \kk$ equation yields after $\tau$-integration
$-\frac {w^2}{2\bar\iota}\left(d_\xx^2+d_\yy^2\right)\log|\bar\iota|=\lambda C$ for a constant
of integration $C$. But $\frac{w^2}{2\bar\iota}\left(d_\xx^2+d_\yy^2\right)\log|\bar\iota|=
\frac 1{2\bar\iota}\left(d_{\bar\xx}^2+d_{\bar\yy}^2\right)\log|\bar\iota|$, so that
a K\"ahler-Einstein metric is obtained precisely when both \Ref{KE-N} and the ODE \Ref{KE} hold.

This completes the proof of Theorem~\ref{ke}.

\subsection{Examples}
First, if $\bar\iota$ is constant, Equation \Ref{KE-N} is satisfied with $C=0$.
The following choices of $f$ and $w$ satisfy the ODE \Ref{KE} with $C=0$, and the inequalities
\Ref{ineq} for some region of $\ta$ values, and thus yield a K\"ahler-Einstein metric given an appropriate Riemannian $3$-manifold $N$, satisfying \Ref{vfield} and \Ref{N-cond} with a specified value
$\al$ and $\bar\iota$ negative. Examples of such $3$-manifolds appearing in \cite[Sec.~8]{a-m} are the
$3$-sphere with the standard metric and $\mathbb{R}^3$ with the truncated pp-wave metric, both having a
global frame satisfying \Ref{N-cond}.

\begin{itemize}
\item{{\bf Vanishing} $\mathbf{\boldsymbol\al}$.}
If $\al=0$, choose $f=1$ and $w(\ta)=(3(a_1p(\ta)+a_2))^{1/3}$, where
\[
p(\ta)=\begin{cases}
\frac{e^{-\lambda \ta}}{-\lambda}\quad \text{if $\lambda\ne 0$}\\
\ta\quad\ \ \ \ \,   \text{if $\lambda=0$}.
\end{cases}
\]
The choice of constants $a_1$, $a_2$ is dictated by the requirements $a_1\ne 0$,
$a_1p(\ta)+a_2>0$ (so $w>0$) and $a_1p'(\ta)>0$ (so $(fw)'>0$).
If $\lambda<0$ these can be satisfied for all values of $\ta$,
but for $\lambda=0$ or $\lambda>0$ only in a subinterval of $\ta$-values. By computing curvatures
of $\gK$ one can show that it is not, in general, of constant sectional curvature.
In the next subsection we give a result on completeness for one of these metrics.
\item{{\bf Negative} $\mathbf{\boldsymbol\al}$.}
If  $\al<0$ and $\lambda=0$, choose $f(\ta)=\ta^{-(1+\al/2)}$ and $w(\ta)=\ta$, limited to the range $\ta>0$. A computation of sectional curvatures indicates these
Ricci flat K\"ahler metrics are in fact flat.
\item{$\mathbf{\boldsymbol\al=-2}$.}
If $\al=-2$ and $\lambda=0$,  choose $f=1$, $w(\ta)=-\tan(x(\ta))$ where $x(\ta)$ solves $x(\ta)=\ta+\tan(x(\ta))$.

We describe the third example in more detail, as it yields a Ricci flat K\"ahler
metric on $S^3\times I$ for an open interval $I$. First, by the implicit function theorem
the zero level set of $h(\ta,x)=\ta+\tan x-x$ is given locally as a function $\ta\to x(\ta)$ near
points $(\ta,x)$ for which $x\ne 2\pi k$. It is easy to calculate that
$(\ta_0,x_0)=(1-\pi/4,-\pi/4)$ lies in this level set and $x'(\ta)=-\cot^2(x(\ta))$ wherever $x(\ta)$ is defined and $\ta\ne k\pi$.

With this one checks that $f=1$, $w(\ta)=-\tan(x(\ta))$ solve the ODE \Ref{KE}, and
near $(\ta_0,x_0)$, $w(\ta)>0$. Furthermore $w'(\ta)=-\sec^2(x(\ta))x'(\ta)=
\sec^2(x(\ta))\cot^2(x(\ta))=\csc^2(x(\ta))> 0$ so that near $(\ta_0,x_0)$, $(fw)'>0$, and the Ricci flat K\"ahler metric is defined as stated on $S^3\times I$, for an appropriate interval $I$ near $\ta_0$.

This metric is not flat: applying our covariant derivative formulas one computes
that
\[
K_{\scriptscriptstyle K}(\xx,\yy)=
-\fr 1f\left[\left(\fr{f'}f+\fr{w'}w\right)\left(\fr{w'}{(fw)'}+1\right)+\fr\al{w}\right],
\]
which for $f=1$ and $\al=-2$ becomes
\[
K_{\scriptscriptstyle K}(\xx,\yy)=
\fr 2w(w'-1),
\]
which is clearly nonzero near $\ta_0$ for $w'$ as above.

\end{itemize}

If $\bar\iota$ is nonconstant, equation~\Ref{KE-N} has, of course, solutions (see
the truncated wave example in subsection~\ref{GR}) for various values of $C$. Then in the ODE~\Ref{KE} one can
combine the term involving $C$ with that containing $\alpha$, and then solve just as
in the case where $C=0$.

\subsection{Completeness}

As shown at the end of Section~8.5 in \cite{a-m}, the induced K\"ahler metric
of any admissible Lorentzian warped product can be written in the form $\gK=ds^2+g_s$,
where $s$ is a certain function of $\ta$, namely $\int\!\!\sqrt{c/2}d\ta$, and $g_s$ is
a metric on $N$. As mentioned there, if $N$ is compact such manifolds are complete
whenever $(\inf s, \sup s)=\mathbb{R}$.

The metrics $g_s$, still written via the variable $\ta$, have the form
$2c(w\hat{\bar{\kk}})^2-fw^{-1}\bar\iota g\big|_\HH=2c(w\hat{\bar{\kk}})^2-
fw\bar\iota\bar g\big|_\HH$ (with the usual meaning of a hatted quantity).
For the particular case of the metrics of Theorem~\ref{ke}, recall that $\bar{k}$ is
in fact Killing, and consider the metric on the two-dimensional quotient space
that pulls back to $\bar\iota \bar g\big|_\HH$. Since the projection of the Lie bracket
$[\bar\xx,\bar\yy]$ to $\HH$ vanishes, a direct computation shows that
equation~\Ref{KE-N} is simply the requirement that the Gauss curvature of this quotient metric is
constant. In this way $\gK$ fits into a well-known ansatz \cite{bb} on line bundles over a Riemann
surface equipped with a metric of constant Gaussian curvature, where equation \Ref{KE} represents the
K\"ahler-Einstein requirement on the line bundle.

For such metrics completeness is well-studied, so we will mention just one case.
Let $(N,\bar{g})$ be a {\em compact} Riemannian $3$-manifold
with a unit length vector field $\bar\kk$ satisfying the assumptions \Ref{vfield}
in the beginning of this section (geodesic, shear-free, with constant twist function
$\bar{\iota}<0$). Assume also conditions \Ref{N-cond} hold with $\al=0$,
so that the universal cover of $N$ is $\mathbb{R}^3$.
As a special case of our first class of examples, set $f(\tau)=1$
and $w(\tau)=-(3e^{-\lambda\ta}/\lambda)^{1/3}$ for $\lambda<0$ constant.
Then by Theorem~\ref{ke} $N\times\mathbb{R}$
admits a K\"ahler-Einstein metric $\gK$. We show first that it does not,
in fact, have constant sectional curvature: calculating with the covariant
derivative formulas in subsection \ref{basic} with $c=(fw)'/w=w'/w$ one sees, for example
that
\[
K_{\scriptscriptstyle K}(\kk,\tT)=-(c'/c)'/c-3c'/c-2c,\qquad
K_{\scriptscriptstyle K}(\xx,\kk)=-c'/(2c)-c/2.
\]
Now $c=-\lambda/3$ so that $K_{\scriptscriptstyle K}(\kk,\tT)=2\lambda/3$
and $K_{\scriptscriptstyle K}(\xx,\kk)=\lambda/6$. In other words these, and
in fact, all frame plane distributions have constant sectional curvature
along the manifold, but comparing, say, the $\kk$, $\tT$ frame plane to
the $\xx$, $\kk$ one at a given point, their sectional curvature values differ.

On the other hand the integral over $\mathbb{R}$ of $\sqrt{c/2}$ is clearly infinite,
so by the result of \cite[subsection~8.5]{a-m}, the metric $\gK$ is complete.

\appendix
\section{Scalar-flat K\"ahler surfaces arising from pp-waves}

We address here the question of whether it is possible to generate scalar flat-K\"ahler
surfaces with symmetry from an ansatz similar to the one used for admissible metrics,
and give examples where the background Lorentzian metric is a pp-wave.

Recall LeBrun's general ansatz \cite{leb} for such K\"ahler metrics,
\[ \text{$\gK=e^uw(dp^2+dq^2)+wdz^2+w^{-1}\theta^2$,} \]
where $p$, $q$, $z$ form a coordinate system for a region in $\mathbb{R}^3$,
and the manifold $M$ is the total space of a circle bundle over such region
(provided the de Rham class of a certain curvature $2$-form associated to the connection $1$-form $\theta$ is integral). Note that $z$ is also a hamiltonian for a holomorphic Killing field,
and the following PDEs hold for $u$, $w>0$:
\be\lb{sfK}
\begin{aligned}
&u_{pp}+u_{qq}+(e^u)_{zz}=0\\
&w_{pp}+w_{qq}+(we^u)_{zz}=0.
\end{aligned}
\end{equation}

With this ansatz in mind, let $\{p,q,z,t\}$ be local coordinates on a manifold $M$,
and fix two smooth functions $u(p,q,z)$, $w(p,q,z)>0$. Consider the 2-form
\begin{equation}\lb{symplt}
\om(a,b):=d(e^\ta \pp^\flat)(a,b)
\end{equation}
where $\tau$ is a smooth function, $\pp$ is a vector field with $1$-form $\pp^\flat$ dual to it with respect to some given semi-Riemannian metric. Let $a$, $b$ take values in frame fields $\kk$, $\tT$, $\xx=\partial_p$, $\yy=\partial_q$ residing in the coordinate neighborhood. Mimicking LeBrun's ansatz we require $d(e^\tau\pp^\flat)(a,b)$ to have the value $e^uw$ when $a=\xx$, $b=\yy$, the value $w$ when $a=\kk$, $b=\tT$, and the value zero for all other pairs
$a$, $b$ with $a\in\{\kk,\tT\}$, $b\in\{\xx,\yy\}$.

Furthermore, we define an almost complex structure by $J\kk=\tT$, $J\xx=\yy$. Note that
the above values on the frame imply that $\om$ is $J$-invariant and symplectic.
We further require
\begin{itemize}
\item $J$ is integrable;
\item $\kk$ preserves $J$;
\item $\kk$ is a hamiltonian vector field for $\om$, with hamiltonian $z$;
\item The following PDEs hold for $u$, $w$:
\begin{align}\lb{laplac}
&d^2_\xx u +d^2_\yy u+(e^u)_{zz}=0,\nonumber\\
&d^2_\xx w +d^2_\yy w+(we^u)_{zz}=0.
\end{align}
\end{itemize}
If these conditions hold, it will follow from \cite{leb} that $\gK=\om(\cdot, J\cdot)$ is a
scalar-flat K\"ahler metric defined in the coordinate domain.

We now demonstrate by verifying these conditions, that examples of this construction hold,
in which the semi-Riemannian background metric is a {\em pp-wave}, given by
\begin{equation}\lb{pp}
g_L=H(x,y,u)du^2+2du\odot dv+dx^2+dy^2.
\end{equation}
We choose our frame as follows:  $\xx=\partial_x+\partial_y+\partial_v$,
$\yy=-\partial_x/2-\partial_y/2+\partial_u-H\partial_v/2$,
$\kk=\partial_v$, $\tT=-3\partial_x-\partial_y-2\partial_v$.
Next we take $\pp=\partial_x-\partial_u+(H+1)\partial_v/2$, so that
\[ \pp^\flat=dx+\frac {1-H}2du-dv.\]

Just as in Remark~\ref{Nintegra}, integrability of $J$ is checked by
verifying $N(\kk,\xx)=0$, where $N$ is the Nijenhuis tensor.
This holds if and only if $3H_x+H_y=0$.
We choose either $H(x,y,u)=e^{x-3y}$ or $H(x,y,u)=x-3y$. In the
latter case $g_L$ is a flat pp-wave and the frame Lie brackets
satisfy the Lie algebra relations of $\mathfrak{nil}_3\times\mathbb{R}$.
We check the remaining conditions only for the latter case,
as the former case is similar.

One calculates that for $\ta=\ta(x,y)$, to obtain the above values of $\om$
on our frame, and hence its $J$-invariance, $\ta_x+\ta_y=0$ must be required
(specifically to have, say, $\om(\kk,\xx)=0$), and we specialize to the case
$\ta=\log\,(y-x)$.

The nonzero values of the of the 2-form on our frame
fields are
\[
\om(\xx,\yy)=e^\ta,\qquad \om(\kk,\tT)=2,
\]
so that we take $w=2$ and $u=\ta-\log 2=\log[(y-x)/2]$
(for $H=e^{x-3y}$ one has instead $\om(\xx,\yy)=e^\ta e^{x-3y})$.
Also, $\iota_\kk\om=\iota_{\partial_v}\om=dy-dx$ and thus the hamiltonian
associated to $\kk$ is $z=y-x$.   Now $e^u=(y-x)/2=z/2$,
so $(e^u)_{zz}=0$, while $(d^2_\xx+d^2_\yy)u=-1/(y-x)^2+1/(y-x)^2=0$.
Thus equations~\Ref{laplac} clearly hold, and $\gK$ is scalar-flat
K\"ahler on the region in $\mathbb{R}^4$ given
by $\{y>x\}$. In fact, for both choices of $H(x,y,u)$, the metric
$\gK$ is also conformally Einstein: $e^{-2\ta}\gK$ is an Einstein metric
with scalar curvature $-12$.

\end{document}